\documentclass[sn-mathphys-num]{sn-jnl}

\usepackage{amssymb}
\usepackage{dsfont}
\usepackage{adjustbox}
\usepackage[symbol]{footmisc}

\newcommand*\dif{\mathop{}\!\mathrm{d}}
\DeclareMathSymbol{\sm}{\mathbin}{AMSa}{"39}
\usepackage{graphicx} 
\DeclareMathSymbol{\smin}{\mathbin}{AMSa}{"39}

\usepackage[dvipsnames]{xcolor}
\definecolor{myRed}{cmyk}{0, 0.7808, 0.7807, 0.0}
\definecolor{myRed2}{cmyk}{0, 0.7808, 0.7807, .2}
\definecolor{myGreen}{cmyk}{.8,0., 1., 0.}
\definecolor{myGreen2}{cmyk}{.8,0., 1., .8}
\definecolor{myBlue}{cmyk}{1.,0., 0., 0.3}
\usepackage{mathtools}
\usepackage[ruled,vlined,linesnumbered]{algorithm2e}
\usepackage{float}
\usepackage{subcaption}	
\usepackage{verbatim}
\usepackage{mathrsfs}
\usepackage{pgfplots}
\usepackage{dsfont}
\usepgfplotslibrary{fillbetween}

\usepackage{amsfonts,amsmath}

\newcommand{\R}{\mathbb{R}}
\newcommand{\C}{\mathbb{F}}

\usepackage{stackengine}
\def\subrangle#1{\stackengine{5pt}{}{$\!\scriptstyle #1$}{U}{l}{F}{F}{L}}
\makeatletter
\let\save@rangle\rangle
\def\rangle{\save@rangle\@ifnextchar_{\expandafter\subrangle\@gobble}{}}
\makeatother
\usepackage{nicefrac, xfrac}

\usepackage{amsthm}
\newtheorem{theorem}{Theorem}[section]

\newtheorem{proposition}[theorem]{Proposition}

\newtheorem{remark}[theorem]{Remark}

\usepackage{array,calc}
\usepackage{tkz-euclide}
\usepackage{standalone}
\usepackage{float}
\usepackage{tikz,pgfplots}
\usepackage{tikz-cd}
\usetikzlibrary{decorations,patterns}
\usetikzlibrary{shadings}
\usetikzlibrary{fit}
\usetikzlibrary{matrix}
\usepackage{filecontents}

\title{Applications of QR-based Vector-Valued Rational Approximation}
\author{
  Simon Dirckx\footnote{Oden institute, University of Texas, \url{simon.dirckx@austin.utexas.edu}},
   Joar Bagge\footnote{KTH Royal Institute of Technology, \url{joarb@kth.se}}, 
   Daan Huybrechs\footnote{Numerical Analysis and Applied Mathematics research unit (NUMA), KU Leuven, \url{daan.huybrechs@kuleuven.be}}
  }
\date{November 2023}
\abstract{
Several applications of the QR-AAA algorithm, a greedy scheme for vector-valued rational approximation, are presented. The focus is on demonstrating the flexibility and practical effectiveness of QR-AAA in a variety of computational settings, including Stokes flow computation, multivariate rational approximation, function extension, the development of novel quadrature methods and near-field approximation in the boundary element method.}
\pgfplotsset{compat=1.18}

\begin{document}
\maketitle
\renewcommand*{\thefootnote}{\arabic{footnote}}
\section{Introduction}
The QR-AAA algorithm, from \cite{PQRAAA}, is a recent development in vector-valued rational approximation by means of greedy optimization. By vector-valued rational approximation, we mean the approximation of a vector-valued function $\mathbf{f}:\mathbb{F}\mapsto \mathbb{F}^n:z\mapsto [f_1(z),\ldots,f_n(z)]$, by means of a vector-valued rational function $\mathbf{r}:\mathbb{F}\mapsto \mathbb{F}^n$. Here $\mathbb{F}=\mathbb{R}$ or $\mathbb{F}=\mathbb{C}$. 

In the ``classical'' SV-AAA scheme from \cite{KarlSVAAA} and \cite{HochmanFastAAA}, a vector-valued rational approximation to $\mathbf{f}$ is built up from samples $[\mathbf{f}_1,\ldots,\mathbf{f}_n] = [f_1(Z),\ldots,f_n(Z)]\in\mathbb{F}^{|Z|\times n}$ of the component functions on a shared discrete sample set $Z$. In QR-AAA one builds a vector-valued rational approximant by first computing an approximate basis $\mathbf{Q}=[\mathbf{q}_1,\ldots,\mathbf{q}_{k}]\in\mathbb{F}^{|Z|\times k}$ of $[\mathbf{f}_1,\ldots,\mathbf{f}_n]$, after which the usual SV-AAA method is applied to (a diagonal scaling of) $\mathbf{Q}$. From the rational approximant to $\mathbf{Q}$, a rational approximant to $\mathbf{f}$ can easily be deduced.

Since in typical applications $k\ll n$, this presents a significantly more efficient method for vector-valued rational approximation. In practice QR-AAA is much faster than SV-AAA, while sacrificing none of the accuracy and robustness. 

Many variants to the basic vector-valued rational approximation schemes from~\cite{KarlSVAAA} and~\cite{HochmanFastAAA} have been proposed in recent years. Some examples include: tangential AAA (tAAA), Block-AAA, stable multi-input AAA (smiAAA), stabAAA, robust surrogate AAA and sketchAAA (\cite{tAAA}, ~\cite{blockAAA}, ~\cite{smiAAA},~\cite{stabAAA},~\cite{robustRatGuttel} and~\cite{sketchAAA}). In contrast to many of these approaches, QR-AAA does not assume or impose additional structure of the vector-valued function nor does it utilize randomized sketching. This significantly simplifies the theoretical derivations and the error estimates, both a priori and a posteriori. In practice though, all these methods have advantages and disadvantages. It depends on the application which, if any, is to be preferred.

In this manuscript, several applications of the QR-AAA framework are presented, which fall into two broad categories:
\begin{enumerate}
    \item Established applications of AAA/SV-AAA that are significantly accelerated through the use of QR-AAA.
    \item New applications that illustrate the flexibility and effectiveness of the QR-AAA approach.
\end{enumerate}
In Section~\ref{sec:qr-aaa}, adapted from~\cite{PQRAAA}, the QR-AAA method is introduced. Afterwards, its effectiveness is demonstrated in Section~\ref{sec:applications} through five applications.

We begin in Section~\ref{sec:Hedgehog} by accelerating the AAA-based \emph{Hedgehog method} of \cite{baggeHedgehog}. The Hedgehog method solves the Stokes equations via a boundary integral equation formulation, in which the dominant computational cost arises from the evaluation of the near-surface Stokes double-layer potential. QR-AAA is used to significantly accelerate this evaluation, leading to improved overall performance of the method.

In Section~\ref{sec:quad}, we introduce a new quadrature scheme for pre-selected classes of functions through vector-valued rational approximation. The main realization here is this: since vector-valued rational approximants share support nodes, SV-AAA can be leveraged to obtain a quadrature scheme using these as quadrature nodes.

In Section~\ref{sec:multivariate}, a partial solution to the problem of practical multivariate rational approximation is presented. We show that for many types of functions, especially smooth functions, a multivariate rational approximation of moderate degree can be obtained using a \textit{Tucker tensor decomposition} generalization of the QR-AAA approach. Additionally, a two-step method is introduced that merges these ideas with the parametric AAA method from~\cite{pAAA}.

In Section~\ref{sec:extension}, multivariate analytic function extension is considered. Using (multivariate) QR-AAA we construct functions that can be evaluated far outside their sample domain, while agreeing to high accuracy with given function samples.

In Section~\ref{sec:BEM}, the parallel version of QR-AAA from \cite{PQRAAA} is applied to the large-scale computation of the wavenumber-dependent near-field of a boundary integral operator for the Helmholtz equation.

It is important to keep in mind the scope of this manuscript. We do not claim that for each of these applications the QR-AAA approach is the only or ``best'' approach. Rather, this collection of applications is meant to demonstrate the flexibility of QR-AAA and inspire the reader to add it to their tool box for the future.

\section{The QR-AAA method}\label{sec:qr-aaa}
The AAA method is a greedy adaptive (scalar) rational approximation scheme. It is practical, easy-to-use, and remarkably robust in that it is reliable for user-defined precisions up to at least $10^{{-8}}$.

Given a sufficiently fine discretization $Z = \{z_1,\ldots,z_N\}\subset \mathbb{F}$ for the domain of a function $f$, the AAA algorithm greedily constructs a rational interpolant in \emph{barycentric form}. Explicitly, the degree $(m-1,m-1)$ rational approximation to $f$ obtained after $m$ iterations of the AAA algorithm is given in the form
$$f(z) \approx r_m(z)=n_{m}(z)\Big/d_m(z)=\left(\sum_{\nu=1}^m\frac{w_{\nu}f(\zeta_{\nu})}{z-\zeta_{\nu}}\middle/\sum_{\nu=1}^m\frac{w_{\nu}}{z-\zeta_{\nu}}\right)$$
with pairwise distinct support points $Z_m=\{\zeta_1,\ldots,\zeta_m\}\subset\C$ and nonzero weights $\{w_{\nu}\}_{\nu}$, subject to $\sum_\nu |w_{\nu}|^2=1$. The expression for $r_m$ is understood to be taken in the `l'H\^opital limit' as $z\to Z_m$.

The greedy nature of this procedure stems from the fact that AAA computes the rational approximant at iteration $m$ by ensuring that
\begin{enumerate}
\item $\zeta_{m}:=\arg\max_{Z\backslash Z_{m-1}} \left|f(z)-r_{m-1}(z)\right|$, and
\item the weights minimize $\|d_{m}(Z \backslash Z_{m})f(Z \backslash Z_{m})-n_{m}(Z \backslash Z_{m})\|_{2}$, the sampled linearized residual in $Z\backslash Z_m$.
\end{enumerate}
Here, it is imposed that $\sum_\nu |w_{\nu}|^2=1$. Setting $Z\backslash Z_{m}=\{\zeta_{m+1},\ldots,\zeta_{N}\}$, this is achieved by computing the compact singular value decomposition (SVD) of the Loewner matrix
$$\mathbf{L}^{(m)}:=\begin{bmatrix}
\frac{f(\zeta_{m+1})-f(\zeta_1)}{\zeta_{m+1}-\zeta_1}&\cdots&\frac{f(\zeta_{m+1})-f(\zeta_m)}{\zeta_{m+1}-\zeta_m}\\
\vdots&\ddots&\vdots\\
\frac{f(\zeta_{|Z|})-f(\zeta_1)}{\zeta_{|Z|}-\zeta_1}&\cdots&\frac{f(\zeta_{|Z|})-f(\zeta_m)}{\zeta_{|Z|}-\zeta_m}
\end{bmatrix}\in\C^{(N-m)\times m}$$
and setting the weight vector $\mathbf{w}:=[w_1,\ldots,w_m]^T$ to be a right singular vector associated with the smallest (not necessarily unique) singular value of $L^{(m)}$. The AAA algorithm provides a remarkably fast, flexible and stable framework for rational approximation (see, e.g., \cite{AAAequi}, \cite{robustRatGuttel}, \cite{KarlSVAAA} and \cite{GuideSaadRatApprox}).

\textit{Vector-valued rational approximation} (also called \textit{set-valued rational approximation}) concerns the approximation of vector-valued functions
$$\mathbf{f}:\mathbb{F}\to\mathbb{F}^n:z\mapsto \mathbf{f}(z)$$
by a vector-valued rational function of type $(m-1,m-1)$ given in barycentric form
\begin{equation}\label{eq:sv-aaa-approx}
\forall z\in \mathbb{F}: \mathbf{f}(z)\approx \mathbf{r}_m(z)=\mathbf{n}_{m}(z)\Big/\mathbf{d}_m(z)=\left(\sum_{\nu=1}^m\frac{w_{\nu}\mathbf{f}(\zeta_{\nu})}{z-\zeta_{\nu}}\middle/\sum_{\nu=1}^m\frac{w_{\nu}}{z-\zeta_{\nu}}\right),  
\end{equation}
as in \cite{KarlSVAAA}. As before we denote by $Z_m=\{\zeta_1,\ldots,\zeta_m\}\subset Z$ the set of chosen support points. Equation~\eqref{eq:sv-aaa-approx} implies that the final degree $(m-1,m-1)$ rational approximations $r_{m,1},\ldots,r_{m,n}$ to $f_1,\ldots,f_n$ share support points and weights, and hence they also share poles.

As in the case of regular AAA, the SV-AAA algorithm from \cite{KarlSVAAA} takes as input samples from $\mathbf{f}$ on a set $Z:=\{z_1,\ldots,z_{N}\}\subset \C$, and greedily selects at every iteration $m$ a new support point $\zeta_m$ by
$$\zeta_{m}:=\arg\max_{Z\backslash Z_{m-1}} \|\mathbf f(z)-\mathbf{r}_{m-1}(z)\|_{p}$$
for some $\|\cdot\|_{p}$-norm\footnote{The two natural choices of the $\|\cdot\|_{p}$-norm are $p=2$ and $p=\infty$.}. Then, it determines the corresponding weights that minimize the linearized residual $\|\mathbf{d}_{m}\mathbf{f}-\mathbf{n}_{m}\|_{p}$ over $Z \backslash Z_{m}$ by computing a singular vector associated with the smallest singular value of a large block Loewner matrix $\mathbb{L}^{(m)}\in \mathbb{F}^{n(N-m)\times m}$. Computing the SVD of $\mathbb{L}^{(m)}$ naively would result in an $\mathcal{O}(m^2n(N-m))$ cost at every iteration. However, in \cite{KarlSVAAA}, a strategic updated SVD scheme for the Loewner matrix is used, reducing the cost for the weight update step at each AAA iteration to $\mathcal{O}(m^2n+m^3)$ flops. Since in typical applications $m\ll n$, this can greatly reduce the overall computation time. It can occur, however, that this updated SVD strategy leads to a loss of stability. For the applications in this manuscript, this does not occur, or can be resolved by switching to a simple regularized version of the update scheme.

As $n$ becomes large (i.e., $\mathbf{f}$ is a very large collections of functions), the computational cost of SV-AAA still becomes prohibitively high. This is the main motivation for QR-AAA, introduced in \cite{PQRAAA}. In short, instead of applying SV-AAA to $\mathbf{f}$, it is applied to a set of basis functions that span the components of $\mathbf{f}$. Section~\ref{sec:interp} provides the details of QR-AAA and the connection of vector-valued rational approximation to interpolative decompositions.
\subsection{Vector-valued rational approximation and interpolative decompositions}\label{sec:interp}
We introduce now two matrices, $\mathbf{F}(Z)\in\C^{N\times n}$ and $\mathbf{R}_m(Z)\in\C^{N\times n}$ by collecting the discretized component functions of $\mathbf{f}$ and $\mathbf{r}_m$ as columns in an array. Concretely:
\begin{equation}\label{eq:defFmat}
    (\mathbf{F}(Z))(i,j) := f_j(z_i) \qquad \text{ and } \qquad (\mathbf{R}_m(Z))(i,j) := r_{m,j}(z_i).
\end{equation}
When the set $Z$ is clear from context we will abbreviate $\mathbf{F}:=\mathbf{F}(Z)$ and $\mathbf{R}_m:=\mathbf{R}_m(Z)$.

As shown in \cite{PQRAAA}, an SV-AAA approximation can be interpreted as an approximate row-interpolative decomposition (RID) of the matrix $\mathbf{F}$ from equation~\eqref{eq:defFmat}. For notational convenience we identify $Z_m\subset Z$ with its corresponding index set in $\{1,\ldots,N\}$.
\begin{theorem}\label{thm:obs1}Suppose $\mathbf{f}:\C\to\C^n$ and $\mathbf{r}_m$ is an SV-AAA approximation of $\mathbf{f}$ with support points $Z_m=\{\zeta_1,\ldots,\zeta_m\} \subset Z=\{z_1,\ldots,z_{N}\}$, such that for $\epsilon>0$ 
\begin{equation}\label{eq:eqobs1}
\mathbf{res}_m:=\sup_{z\in Z}\|\mathbf{f}(z)-\mathbf{r}_m(z)\|_{p}<\epsilon.
\end{equation}
Let $\mathbf{F}:=\mathbf{F}(Z)$ and $\mathbf{R}_m:=\mathbf{R}_m(Z)$ be as in equation~\eqref{eq:defFmat}. Then the matrix $\mathbf{R}_m$ is of rank at most $m$, and can be written as
\begin{equation}\label{eq:IDHthm}
\mathbf{R}_m=\mathbf{H}_m \mathbf{F}(Z_m,:)
\end{equation}
where $\mathbf{H}_m\in\mathbb{F}^{N\times m}$ is such that
$$\mathbf{H}_m(i,\nu) = \left.\frac{w_{\nu}}{z_i-\zeta_{\nu}}\middle/\sum_{\nu=1}^m\frac{w_{\nu}}{z_i-\zeta_{\nu}}\right. $$
which as before is evaluated as a limit for $z_{i}\in Z_m$.
The decomposition in equation~\eqref{eq:IDHthm} constitutes an approximate RID of $\mathbf{F}$, with respect to the $\|\cdot\|_{p,\infty}$-norm\footnote{Meaning $\|\mathbf{F}-\mathbf{R}_m\|_{p,\infty}:=\max_{i}\|\mathbf{F}(i,:)-\mathbf{R}_m(i,:)\|_p<\epsilon$}.
\end{theorem}
Note that $\|\mathbf{F}-\mathbf{R}_m\|_{p,\infty}<\epsilon$ implies $\|\mathbf{F}-\mathbf{R}_m\|_{\text{max}}<\epsilon$. The matrix $\mathbf{H}_m$ will be referred to as the \textit{barycentric matrix} associated with $\{\zeta_1,\ldots,\zeta_m\}$ and $\{w_1,\ldots,w_m\}$.

The principle behind QR-AAA is straightforward. Suppose that, with $0<\epsilon_2$ sufficiently small 
$$\mathbf{F}\mathbf{P} =\begin{bmatrix}
    \mathbf{Q}_1&\mathbf{Q}_2
\end{bmatrix}
\begin{bmatrix}
    \mathbf{C}_1\\\mathbf{C}_2
\end{bmatrix}:=\mathbf{Q}_1\mathbf{C}_1+\mathbf{\Delta}_2,\quad \|\mathbf{\Delta}_2\|_{\text{max}}<\mathbf{\epsilon}_2,$$
with $\mathbf{Q}_1\in \mathbb{F}^{n\times k}$, is obtained using a column-pivoted QR-decomposition\footnote{$\mathbf{C}$ is used instead of $\mathbf{R}$ to avoid confusion with $\mathbf{R}_m$. Exactness of the QR decomposition can be relaxed by absorbing its error into $\mathbf{\Delta}_1$.}. If we apply SV-AAA to $\mathbf{Q}_1$, we obtain an approximation $\mathbf{H}_m\mathbf{Q}_1(Z_m,:)$ which satisfies
$$\mathbf{Q}_1 = \mathbf{H}_m\mathbf{Q}_1(Z_m,:)+\mathbf{\Delta}_1,\quad \|\mathbf{\Delta}_1\|_{p,\infty}<\epsilon_1.$$ 
Then it can be shown (using H\"older's inequality) that
$$\|\mathbf{F}-\mathbf{H}_m\mathbf{F}(Z_m,:)\|_{\text{max}}< \epsilon_2+\frac{\epsilon_1k}{\sqrt[p]{k}}\|\mathbf{F}\|_{\text{max}}+\|\mathbf{H}_m\mathbf{\Delta}_2(Z_m,:)\|_{\text{max}}.$$
We assume $\mathbf{H}_m$ is ``well-behaved'', meaning 
$$\|\mathbf{\Delta}_2(Z_m,:)\|_{\text{max}}<\epsilon_2 \quad\text{ implies }\quad \|\mathbf{H}_m\mathbf{\Delta}_2(Z_m,:)\|_{\text{max}}<C\epsilon_2$$ with a modest $C$ (``small sample values are mapped to small functions''). This means it cannot introduce \textit{spurious poles}, i.e., poles inside the domain or isolated poles close to the domain not corresponding to near-singularities in the data. The case of poles clustering towards the domain is much harder to analyze, but in practice the barycentric matrices are still well-behaved in this case.

The above means that applying the SV-AAA algorithm to $\mathbf{Q}_1$ can be used to obtain a vector-valued approximation of $\mathbf{F}$. In practice, SV-AAA is applied to $\mathbf{Q}_1\mathrm{diag}(\mathbf{C}_1)$, to avoid overfitting the trailing, typically oscillatory, columns of $\mathbf{Q}_1$. The QR-AAA approach is summarized in Algorithm~\ref{alg:QRAAA}. We use \texttt{MATLAB} notation where convenient.

\begin{algorithm}
\caption{Vector-valued rational approximation by QR-AAA}\label{alg:QRAAA}
\SetKwInOut{Input}{Input}
\SetKwInOut{Output}{Output}
\SetKw{Init}{init}{}{}
\SetAlgoLined
\DontPrintSemicolon
\Input{Sample points $Z$, matrix $\mathbf{F}(Z)\in \mathbb{F}^{N\times n}$ of the vector-valued function $\mathbf{f}(z)$ sampled on $Z$ (see equation~\ref{eq:defFmat}),  tolerances $\texttt{tol\_qr},\texttt{tolAAA}$}
\Output{Vector-valued rational approximant $\mathbf{r}_m(z)$.}
 $[\mathbf{Q},\mathbf{C},\sim] = \texttt{qr}(\mathbf{F},\texttt{`econ'})$\;
    \texttt{dc} $=$ \texttt{abs(diag(}$\mathbf{C}$\texttt{))}\;
    \texttt{k} $=$ \texttt{sum(dc>tol\_qr*dc(1))}\;
    $\mathbf{Q}_{k} = \mathbf{Q}$\texttt{(:,k).*dc(1:k)'}\;
    \tcp{Compute SV-AAA approx, keep only support points and weights}
    \texttt{[ $\sim$ , $Z_{m}$, $\mathbf{w}_{m}$]} $=$ \texttt{aaa\_sv($\mathbf{Q}_{k}$,$Z$,`tol',tolAAA)}\;
    \tcp{With r as in equation~\ref{eq:sv-aaa-approx}}
    $\mathbf{r}_m$ = \texttt{@(z)r($Z_m$,$\mathbf{w}_{m}$,$\mathbf{f}(Z_m)$)($z$)}
\end{algorithm}
In a performant implementation of the QR-AAA algorithm, only the leading $k\times k$-subblock of $\mathbf{C}$ is computed, using \textit{early-exit} column-pivoted Householder QR. This is the version included in the \texttt{C++} package PQR-AAA, available at \url{https://github.com/SimonDirckx/pqraaa}. For the purposes of this manuscript however, a simple \texttt{MATLAB} version of the QR-AAA algorithm, as well as the code for the experiments in Sections~\ref{sec:Hedgehog}-\ref{sec:extension}, is available at \url{https://github.com/SimonDirckx/QR-AAA-Applications}.

\section{Applications of QR-AAA}\label{sec:applications}
\subsection{A vector-valued Hedgehog method}\label{sec:Hedgehog}
When solving the Stokes equations using boundary integral methods, one often needs to evaluate the \textit{double layer potential} $\mathcal{D}[\mathbf{q}]$, given by
$$(\mathcal{D}_{i}[\mathbf{q}])(\mathbf{x}) = \sum_{j,k=1}^3\int_{\Gamma}T_{ijk}(\mathbf{x},\mathbf{y})q_{j}(\mathbf{y})\,n_{k}(\mathbf{y}) dS(\mathbf{y})$$
with the \textit{stresslet} $T_{ijk}$ defined by 
$$T_{ijk}(\mathbf{x},\mathbf{y}) = -6\frac{(x_i-y_i)(x_j-y_j)(x_k-y_k)}{\|\mathbf{y}-\mathbf{x}\|^5},$$
$i,j,k\in\{1,2,3\}$ and $\Gamma$ the boundary of a domain of interest. Here $\mathbf{q}=[q_1,q_2,q_3]$ is a \textit{layer density} on $\Gamma$. Using $\mathcal{D}[\mathbf{q}]$, the associated flow field $\mathbf{u}$ away from and on the boundary $\Gamma$ can be computed after $\mathbf{q}$ has been obtained. More details can be found in \cite{baggeStokes}.

\begin{figure}
\begin{subfigure}[valign = t]{.55\linewidth}
\centering
\includegraphics[width=\linewidth]{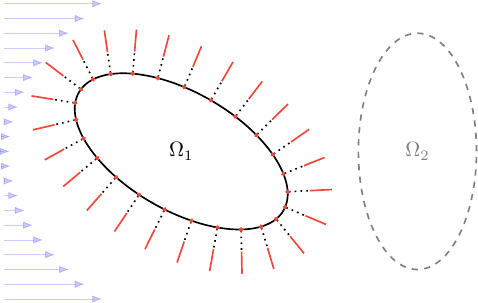}
\caption{Rigid body $\Omega_1$ in Stokes flow with Hedgehog spikes illustrated.}
\end{subfigure}
\hfill
\begin{minipage}[valign=t]{.4\linewidth}
    \textit{Diagram of two rigid bodies in a 3D Stokes flow with paraboloidal background flow (indicated in blue). On $\Gamma_1$ the Hedgehog spikes are shown. Each of these is divided in three parts: a far-field part (red line), a near-field part (dashed) and an on-surface point (red). The double-layer density is fitted on the red parts and interpolated to the near-field.}
\end{minipage}
        \caption{Diagram illustrating the Hedgehog set-up under consideration.}
    \label{fig:Stokes}
\end{figure}

For instance, one can consider the case of rigid bodies in a paraboloid background flow, as depicted in Figure~\ref{fig:Stokes}. In this set-up, the rigid body motion of two floaters, along with their associated layer densities $\mathbf{q}_1=[q_{1,1},q_{1,2},q_{1,3}]$ on $\Gamma_1$ and $\mathbf{q}_2=[q_{2,1},q_{2,2},q_{2,3}]$ on $\Gamma_2$ are computed. Subsequently the flow field contribution $\mathbf{u}_1$ from $\Gamma_1$ is computed using (amongst other quantities) the double layer potential $\mathcal{D}[\mathbf{q}_1]$. In this application we will focus only on the computation of $\mathcal{D}[\mathbf{q}_1]$, but the same method can be used to obtain $\mathcal{D}[\mathbf{q}_2]$.

The main difficulty in computing $\mathcal{D}[\mathbf{q}_1]$ comes from the fact that the stresslet contains a strong singularity for $\mathbf{x}\rightarrow\Gamma_1$, which renders regular quadrature methods insufficient. While for the on-surface evaluation so-called \textit{QBX quadrature} can be very efficient (see \cite{baggeStokes,QBX,QBXclose,KlintebergQBX}), the evaluation of $\mathcal{D}[\mathbf{q}_1]$ in a boundary layer close to the surface can still be prohibitively expensive. Various QBX acceleration schemes to overcome this cost have been proposed (see, e.g., \cite{WalaQBX}), but in this work the focus will be on the method of \textit{line interpolation}, also known as the \textit{Hedgehog method} introduced in \cite{YINGHedgehog}.

In the variant of the Hedgehog method from \cite{baggeHedgehog}, depicted in Figure~\ref{fig:Stokes}, the double layer potential $\mathcal{D}[\mathbf{q}_1]$ is fitted far away from \textbf{and} on the boundary $\Gamma_1$ and subsequently interpolated to a region close to the boundary. Practically, this is done by introducing a collection of \textit{spikes} $\{s_i\}_{i=1}^n$, each of which is an affine mapping (say, $\mathbf{A}_i$) of $\{0\}\cup [\alpha,\beta]$. For each of these spikes, $\mathcal{D}[\mathbf{q}_1]$ is computed on $\mathbf{A}_i([\alpha,\beta])$ using regular upsampled quadrature. The double layer potential at the boundary (i.e., at $\mathbf{A}_i(0)$) is computed using QBX. Finally, AAA rational approximation is used for the fitting of the resulting data on each spike. As the number of spikes grows, so does the overall computational cost. This limits the accuracy that can be achieved in the volume and as such motivates the use of QR-AAA.

To be precise, in our set-up QR-AAA is applied to a 3D Stokes problem involving two rigid floating bodies. We construct an approximation for the collection of spike functions $\mathcal{V}:=\{\mathbf{v}_i\}_{i=1}^n$ with
$$\mathbf{v}_i:\{0\}\cup[\alpha,\beta]\rightarrow \mathbb{R}^3:x\mapsto (\mathcal{D}[\mathbf{q}_1])(\mathbf{A}_i(x))$$
where $\mathbf{A}_i$ is the affine mapping corresponding to the $i$-th Hedgehog spike. Once a discretization $Z$ of $\{0\}\cup [\alpha,\beta]$ is chosen, this corresponds to applying the QR-AAA algorithm to the matrix $\mathbf{V}=[\mathbf{V}_1,\mathbf{V}_2,\mathbf{V}_3]\in \mathbb{R}^{|Z|\times 3n}$ with $\mathbf{V}_d(i,j)=v_{j,d}(z_i)$\footnote{The construction of $\mathbf{V}$ is the dominant cost in our \texttt{MATLAB} code. In the case of Figure~\ref{tab:StokesResults}, the time to construct $\mathbf{V}$ is around 58.8~s for $|Z|=10^2$, and 573.9~s for $|Z|=10^3$. This cost depends on the underlying potential evaluator for $\mathcal{D}[\mathbf{q}_1]$. Our potential evaluator is not optimized, and runs only on a single core. It could be improved by implementing a parallel version in a high-performance programming language (the computation is embarrassingly parallelizable), but this is not the focus of this section. 
}.
In our set-up, $n=512$ was chosen, and $Z$ was set to be equispaced points in $[\alpha,\beta]=[0.1,0.2]$, together with the origin. We report our findings in Figure~\ref{fig:StokesResults} and the table in Figure~\ref{tab:StokesResults}.

\begin{figure}
    \centering
    \begin{subfigure}[b]{.52\linewidth}
    \begin{tikzpicture}
    \begin{axis}[
    width=\linewidth,
    ymode=log,
    log basis y=10,
    xlabel={Distance to $\Gamma$},
    ylabel={$\|\mathbf{v}(x)-\hat{\mathbf{v}}(x)\|_2/\|\mathbf{v}(x)\|$},
    ymin = 1e-10,
    ymajorgrids=true,
    xmajorgrids=true,
    grid style=dashed,
    xtick = {0,.02,.04,.06,.08,.1},
    xticklabels={$0$,$.02$,$.04$,$.06$,$.08$,$.1$},
    legend style={fill=none}
]
\addplot[
    color=black,dashed,forget plot
    ] table [x index = 0 , y index = 1,col sep=comma]{aaaErr.tex};
\addplot[
    color=black,dashed,forget plot
    ] table [x index = 0 , y index = 2,col sep=comma]{aaaErr.tex};
\addplot[
    color=black,dashed,forget plot
    ] table [x index = 0 , y index = 3,col sep=comma]{aaaErr.tex};
\addplot[
    color=black,dashed,forget plot
    ] table [x index = 0 , y index = 4,col sep=comma]{aaaErr.tex};
\addplot[
    color=black,dashed
    ] table [x index = 0 , y index = 5,col sep=comma]{aaaErr.tex};
    \addlegendentry{AAA}
\addplot[
    color=black,forget plot
    ] table [x index = 0 , y index = 1,col sep=comma]{qraaaErr.tex};
\addplot[
    color=black,forget plot
    ] table [x index = 0 , y index = 2,col sep=comma]{qraaaErr.tex};
\addplot[
    color=black,forget plot
    ] table [x index = 0 , y index = 3,col sep=comma]{qraaaErr.tex};
\addplot[
    color=black,forget plot
    ] table [x index = 0 , y index = 4,col sep=comma]{qraaaErr.tex};
\addplot[
    color=black
    ] table [x index = 0 , y index = 5,col sep=comma]{qraaaErr.tex};
    \addlegendentry{QR-AAA}
    \end{axis}
    \end{tikzpicture}
    \caption{Relative error in $[0,\alpha]$ for the $5$ Hedgehog spikes with the largest relative error, for AAA (dashed) and QR-AAA (solid). (Note that the 5 curves for QR-AAA coincide.)}
    \label{fig:StokesResults}
    \end{subfigure}
    \hfill
    \begin{subfigure}[b]{.45\linewidth}
        \begin{tabular}{c|c|c}
    &\small $|Z|=10^2$&\small $|Z|=10^3$\\[5pt]
    \hline
    &&\\
    \small AAA&2.608s&25.637s\\[5pt]
    \small QR-AAA&0.054s&0.303s
    \end{tabular}
    \vspace{2cm}
    \caption{Timings for AAA and QR-AAA for $n=512$ and $|Z|\in\{100,1000\}$.}
    \label{tab:StokesResults}
    \end{subfigure}
    \label{fig:NumericalResultsStokes}
    \caption{Numerical results for QR-AAA as compared to AAA for $512$ Hedgehog spikes, with $[\alpha,\beta]=[.1,.2]$. Left: relative errors for the Hedgehog spikes that have the largest relative interpolation error in $[0,\alpha]$, for $|Z|=100$. Right: timings for the regular Hedgehog method (each spike approximated separately) and for the QR-AAA accelerated version.}
\end{figure}
Let us make a couple of notes:
\begin{enumerate}
    \item Both the AAA and the QR-AAA approximations are limited by the QBX accuracy at the boundary, which is $\sim 10^{-5}$. As we can see, both the AAA approach and the QR-AAA approach also reach approximately $5$ digits of accuracy, with the largest error reached closest to $\Gamma_1$.
    \item QR-AAA was run at a tolerance $10^{-12}$, and the rank of the resulting $\mathbf{Q}_k$ was $10$. At this tolerance, the degree of the vector-valued rational approximant was $8$. For comparison, the maximum degree for the individual AAA approximants was $3$.
    \item Due to the rank of $\mathbf{Q}_k$ being so low, QR-AAA is significantly faster than using regular AAA; in both cases a factor of more than $50$ is achieved. 
\end{enumerate}
In conclusion, QR-AAA achieves a significant speed-up without sacrificing accuracy.

\subsection{Rational quadrature}\label{sec:quad}
Consider the case of generalized Gaussian quadrature (see \cite{GaussQuad} and \cite{QuadDaanRonald}) where we want to construct a set of quadrature \textit{nodes} $\{x_{\nu}\}_{\nu=1}^{m}$ and \textit{weights} $\{c_{\nu}\}_{\nu=1}^{m}$ such that
$$I[f]=:\int_{0}^{1}f(x)\dif x \approx Q[f]=: \sum_{\nu=1}^{m}c_{\nu}f(x_{\nu}).$$
for $f$ in some suitable class of functions. One example of such a class is the functions that are smooth away from zero, but possibly have an algebraic singularity at zero. Examples include $f(x)=x\sqrt{x}$, $f(x)=x^{\pi}$, $f(x) = x^{3/2}\cos(4x)$ and so on. Another interesting class of functions (see \cite{QuadDaanRonald}) are those of the form $f(x) = u(x)+v(x)\log(x)$ where $u$ and $v$ are smooth, or even of the first class, i.e. containing an algebraic singularity at the origin.

We propose the following strategy, inspired by the approach in \cite{GaussQuad} and the principle of direct rational quadrature derived in \cite{BaryQuad}:
\begin{enumerate}
    \item Select a large set of functions $\Psi=\{f_1,\ldots,f_n\}$ that is expected to approximately span the full class of functions. 
    \item Select a sufficiently fine sampling grid $Z\subset [0,1]$ and apply QR-AAA to the sampled functions $[\mathbf{f}_1,\ldots,\mathbf{f}_n]$.
    \item Set the quadrature nodes $\{x_{\nu}\}_{\nu=1}^m$ to be the selected support points $\{\zeta_1,\ldots,\zeta_m\}$.
    \item With the resulting vector-valued rational approximation as in equation~\eqref{eq:sv-aaa-approx}, compute the quadrature weights $\{c_{\nu}\}_{\nu=1}^{m}$ as
    \begin{equation*}
        c_\nu:=\int_0^1\left( \frac{w_{\nu}}{z-\zeta_{\nu}}\middle/\sum_{\nu}\frac{w_{\nu}}{z-\zeta_{\nu}}\right)\dif z
    \end{equation*}
\end{enumerate}
In this way we have that, at least if the functions in $\{f_1,\ldots,f_N\}$ are well-approximated by QR-AAA, their integral satisfies
\begin{align*}
    \int_{0}^{1}f_{i}(x)\dif x &\approx \int_{0}^{1} \left(\sum_{\nu=1}^m\frac{w_{\nu}f_i(\zeta_{\nu})}{z-\zeta_{\nu}}\middle/\sum_{\nu=1}^m\frac{w_{\nu}}{z-\zeta_{\nu}}\right)\dif z\\
    &=\sum_{\nu=1}^mf_i(\zeta_{\nu})\int_0^1\left( \frac{w_{\nu}}{z-\zeta_{\nu}}\middle/ \sum_{\nu}\frac{w_{\nu}}{z-\zeta_{\nu}}\right)\dif z\\
    &=\sum_{\nu=1}^mf_i(x_{\nu})c_{\nu}
\end{align*}
This holds by virtue of the fact that we deliberately sought out an approximation with shared support points.

It should be stressed that this is \textit{not} a Gaussian quadrature rule; we pre-select candidate quadrature nodes and the final quadrature rule selects a quasi-optimal subset from these, with a set of weights derived from vector-valued rational approximation. There has been some research towards computing ``proper'' Gaussian quadrature rules based on rational approximation. In \cite{GAUTSCHI2001111}, three-term recurrence coefficients for known Gauss rules are updated to include rational functions with pre-assigned poles. In \cite{VanDeunQuad}, recurrence coefficients are built up directly from orthogonal rational functions with pre-assigned poles. Unfortunately, both approaches suffer from severe instabilities as the poles approach the interval of interest, making them unsuited for the classes of functions considered in this section. It is an interesting open question whether these instabilities can be overcome. If this is the case, QR-AAA can be used to approximate the poles for a class of functions, after which an associated (approximate) Gaussian quadrature rule can be constructed using either of these approaches. 

In \cite{HorningQuad}, a different strategy is used. There, the problem of quadrature with respect to a given weight function is solved by constructing a rational approximant of the Cauchy transform of the given weight function. The poles of this rational approximant form the quadrature nodes, while the residues form the quadrature weights. However, this construction is limited to integrands that are analytic in a sufficiently large neighborhood surrounding the interval of integration. Additionally, this differs from the approach we take here in that it does not allow one to select a function class of interest for which the quadrature rule is sufficiently accurate.

Our approach is similar to, but distinct from, the approach in~\cite[Section 6.2]{RKfit}. There, an exponential integrator is set up using shared pole approximations for a set of rational functions which approximate time-stepping operators. This is done using rational Krylov based approximations, and the final approximant is tailored specifically to the time-stepping problem. By contrast, we employ QR-AAA and consider the general quadrature problem for a provided class of functions.

This section reveals a reciprocal relationship with data-driven reduced order modeling: there, reduced models are constructed using numerical quadrature for judiciously chosen contour integrals, see~\cite{doi:10.1137/20M1389303} and~\cite{doi:10.1137/21M1411081}. Here we proceed in the opposite direction: we start from vector-valued rational approximation and construct a quadrature rule.

\subsubsection{Computing the quadrature weights}\label{sec:weightComp}
After the support points $Z_m$ have been computed using QR-AAA, the computation of the quadrature weights $\{c_{\nu}\}_{\nu=1}^{m}$ defined by
\begin{equation*}
    c_\nu:=\int_0^1\left( \frac{w_{\nu}}{z-\zeta_{\nu}}\middle/\sum_{\nu}\frac{w_{\nu}}{z-\zeta_{\nu}}\right)\dif z
\end{equation*}
can be accomplished in two ways. The first and most direct way is to select a sufficiently accurate conventional quadrature rule. Since we are primarily interested in functions with a singularity at the origin (for which it is known that QR-AAA poles tend to cluster towards the singularity), an adaptive Gauss-Kronrod quadrature rule was used in our experiments.

However, if the exact integrals for the selected set of functions $\Psi = \{f_1,\ldots,f_n\}$ are known, we can leverage this fact for improved accuracy. Let $[I_1,\ldots,I_n]$ denote the array of integrals of $\{f_1,\ldots,f_n\}$ and $\mathbf{c}^T=[c_1,\ldots,,c_m]$. Then, if the QR-AAA approximation is sufficiently accurate, we have 
\begin{equation*}
    \mathbf{c}^T\mathbf{F}(Z_m,:)\approx [I_1,\ldots,I_n].
\end{equation*}
This means we have an alternate way of computing $\mathbf{c}$:
\begin{equation*}
    \mathbf{c}^T := [I_1,\ldots,I_n]\mathbf{F}(Z_m,:)^{\dagger}.
\end{equation*}
As will be shown in Section~\ref{sec:numExpQuad}, this method of computing the quadrature weights is often slightly more reliable than the Gauss-Kronrod based version. However, one should keep in mind that this technique can only be used if the exact integrals of the functions in $\Psi$ are known.
\subsubsection{Numerical experiments}\label{sec:numExpQuad}
\begin{figure}
    \centering
    \begin{minipage}[b]{.65\linewidth}
    \begin{tikzpicture}
    \begin{axis}[
    width=\linewidth,
    xmode=log,
    ymode=log,
    log basis y=10,
    log basis x=10,
    xlabel={$\mathbf{tol}$},
    ylabel={$|\,I[f]-Q[f]\,|$},
    ymajorgrids=true,
    xmajorgrids=true,
    grid style=dashed,
    legend style={fill=none},
    legend pos={north west}
]
\addplot[
    color=black
    ] table [x index = 0 , y index = 1,col sep=comma]{quad_errors.tex};
    \addlegendentry{$x^{\alpha}\sin(x)$}
\addplot[
    color=black,mark=square*
    ] table [x index = 0 , y index = 2,col sep=comma]{quad_errors.tex};
    \addlegendentry{$x^{\alpha}\cos(x)$}
\addplot[
    color=black,mark=triangle*
    ] table [x index = 0 , y index = 3,col sep=comma]{quad_errors.tex};
    \addlegendentry{$J_0(\alpha x)$}
\addplot[
    color=myRed,dashed,forget plot
    ] table [x index = 0 , y index = 1,col sep=comma]{quad_errors_use_exact_int.tex};
\addplot[
    color=myRed,dashed,mark=square*,forget plot
    ] table [x index = 0 , y index = 2,col sep=comma]{quad_errors_use_exact_int.tex};
\addplot[
    color=myRed,dashed,mark=triangle*,forget plot
    ] table [x index = 0 , y index = 3,col sep=comma]{quad_errors_use_exact_int.tex};
    \end{axis}
    \end{tikzpicture}
    \end{minipage}
    \hfill
    \begin{minipage}[b]{.32\linewidth}
        \textit{For each of the selected families of functions, the range $\alpha\in [0,10]$ was chosen. The maximum error over $\alpha$ at the requested tolerance is reported (this occurs typically at the largest $\alpha$). The quadrature order increases linearly from $25$ to $44$.}
        \vspace{5em}
    \end{minipage}
    \caption{Quadrature error for the QR-AAA based quadrature (class $\Psi_1$) over the requested tolerance $\mathbf{tol}$ for various function classes. Quadrature errors with weights computed using exact integrals (the second method from Section~\ref{sec:weightComp}) are shown in red.}
    \label{fig:quadResults}
    \end{figure}

\begin{figure}
    \centering
    \begin{subfigure}[t]{.49\linewidth}
    \begin{tikzpicture}
\tikzset{mark options={mark size=1}}
\begin{axis}[%
width=\linewidth]
\addplot[mark = *,only marks]
table[x index = 0,y expr =\thisrowno{0}*0,col sep=comma] {quad_supp_exact_int.tex};
\end{axis}
\end{tikzpicture}
\caption{Quadrature nodes for the QR-AAA approximation of $\Psi_2$}
\label{subfig:polesquad1}
\end{subfigure}
\hfill
\begin{subfigure}[t]{.49\linewidth}
    \begin{tikzpicture}
\tikzset{mark options={mark size=1}}
\begin{axis}[%
width=\linewidth]
\addplot[mark = *,only marks]
table[x index = 0,y expr =\thisrowno{0}*0,col sep=comma] {quad_sing_supp_exact_int.tex};
\end{axis}
\end{tikzpicture}
\caption{Quadrature nodes for the QR-AAA approximation of $\Psi_2$}
\label{subfig:polesquad2}
\end{subfigure}
    \caption{Scatter plot of the QR-AAA poles selected for the approximation of $\Psi_1$ and $\Psi_2$, with the QR-AAA tolerance set to $10^{-10}$.}
    \label{fig:poles}
\end{figure}
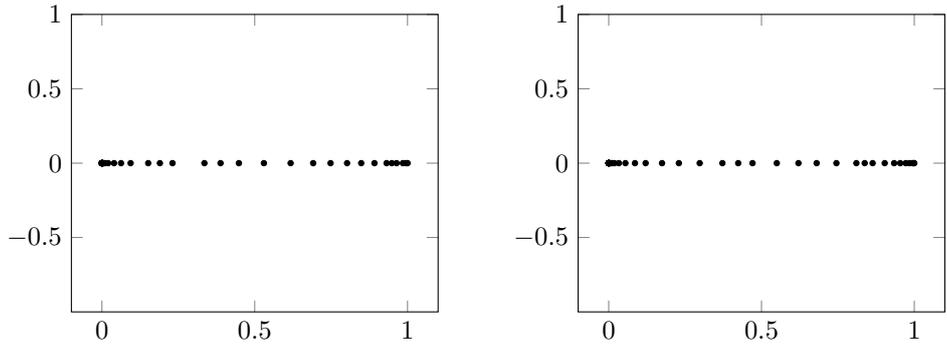

We ran our QR-AAA based quadrature scheme on two families of functions:
$$\Psi_1:=\{x^{\alpha}|\,\alpha\in [0,50]\}$$
and
$$\Psi_2:=\{x^{\alpha}|\,\alpha\in [0,50]\}\cup \{x^{\beta}\log(x)|\,\beta\in [0,1/2)\}$$
Since, for $\Psi_2$, we are interested in the singular behavior around $\beta=0$, we chose a log-spaced discretization of for the parameter $\beta$. For $\alpha$, in both cases, a fine equispaced discretization ($n_{\alpha}=5000$) was chosen. For $\Psi_1$, where there is no singularity at the origin, we chose $Z$ to be a sufficiently fine ($|Z|=5000$) second-kind Chebyshev grid on $[0,1]$. For $\Psi_2$, where there is a singularity at the origin, we chose $Z$ to be a sufficiently fine ($|Z|=5000$) grid of Legendre points. In both cases, the grid was augmented with $50$ log-spaced points clustering towards a point close to the origin.\footnote{Clustering completely towards the origin results in instabilities due to the strong clustering of the resulting poles. In our case we chose to cluster towards $10^{-8}$. It is an interesting open question what the optimal strategy is.} The accuracy of QR-AAA based quadrature with these set-ups, run at various tolerances for the vector-valued rational approximation, are reported for some interesting classes of test functions in Figure~\ref{fig:quadResults} and Figure~\ref{fig:quadSingResults}. Here $J_0(x)$ and $Y_0(x)$ are the order zero Bessel functions of the first and second kind respectively. The resulting quadrature nodes of the QR-AAA approximations, at the requested tolerance $10^{-10}$, are shown in Figure~\ref{fig:poles}.

In all cases, the QR-AAA based quadratures reach an accuracy that agrees well with the requested tolerance. Typically, using the analytically known integrals for the functions in $\Psi_1$ and $\Psi_2$ results in a more accurate quadrature rule. This is especially true for $\Psi_2$, where the functions to be integrated can contain a logarithmic singularity at the origin. We set the tolerance for the QR decomposition to be $10^{-13}$ in all cases. For $\Psi_1$, the rank obtained by the QR-decomposition was $37$. For $\Psi_2$ it was $40$.

\begin{figure}
    \centering
    \begin{minipage}[b]{.65\linewidth}
    \begin{tikzpicture}
    \begin{axis}[
    width=\linewidth,
    xmode=log,
    ymode=log,
    log basis y=10,
    log basis x=10,
    xlabel={$\mathbf{tol}$},
    ylabel={$|\,I[f]-Q[f]\,|$},
    ymajorgrids=true,
    xmajorgrids=true,
    grid style=dashed,
    legend style={fill=none},
    legend pos={south east}
]
\addplot[
    color=black,mark=*
    ] table [x index = 0 , y index = 4,col sep=comma]{quad_sing_errors.tex};
    \addlegendentry{$\cos(\alpha x)\log(x)$}
\addplot[
    color=black
    ] table [x index = 0 , y index = 1,col sep=comma]{quad_sing_errors.tex};
    \addlegendentry{$x^{\alpha}\log(x)$}
\addplot[
    color=black,mark=triangle*
    ] table [x index = 0 , y index = 3,col sep=comma]{quad_sing_errors.tex};
    \addlegendentry{$Y_0(\alpha x)$}
\addplot[
    color=myRed,dashed,mark=*,forget plot
    ] table [x index = 0 , y index = 4,col sep=comma]{quad_sing_errors_use_exact_int.tex};
\addplot[
    color=myRed,dashed,forget plot
    ] table [x index = 0 , y index = 1,col sep=comma]{quad_sing_errors_use_exact_int.tex};
\addplot[
    color=myRed,dashed,mark=triangle*,forget plot
    ] table [x index = 0 , y index = 3,col sep=comma]{quad_sing_errors_use_exact_int.tex};
    \end{axis}
    \end{tikzpicture}
    \end{minipage}
    \hfill
    \begin{minipage}[b]{.32\linewidth}
        \textit{For each of the selected families of functions, the range $\alpha\in [0,10]$ was chosen. The maximum error over $\alpha$ at the requested tolerance is reported (this occurs typically at the largest $\alpha$). The quadrature order increases linearly from $26$ to $47$.}
        \vspace{5em}
    \end{minipage}
    \caption{Quadrature error for the QR-AAA based quadrature (class $\Psi_2$) over the requested tolerance $\mathbf{tol}$ for various function classes. Quadrature errors with weights computed using exact integrals (the second method from Section~\ref{sec:weightComp}) are shown in red.}
    \label{fig:quadSingResults}
    \end{figure}

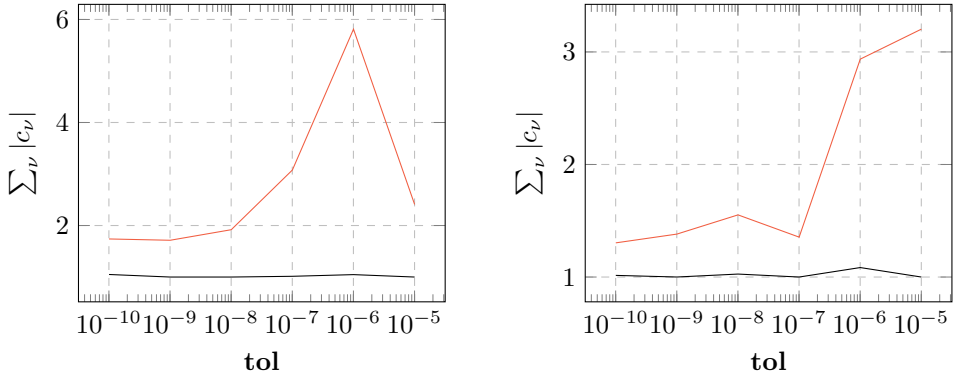
\begin{figure}
    \centering
    \begin{subfigure}[t]{.49\linewidth}
    \begin{tikzpicture}
    \pgfplotsset{
  /pgfplots/xlabel near ticks/.style={
     /pgfplots/every axis x label/.style={
        at={(ticklabel cs:0.5)},anchor=near ticklabel
     }
  },
  /pgfplots/ylabel near ticks/.style={
     /pgfplots/every axis y label/.style={
        at={(ticklabel cs:0.5)},rotate=90,anchor=near ticklabel}
     }
  }
\begin{axis}[%
width=\linewidth,
    xmode=log,
    log basis x=10,
    xlabel={$\mathbf{tol}$},
    ylabel={$\sum_{\nu}|c_{\nu}|$},
    ylabel near ticks,
    ymajorgrids=true,
    xmajorgrids=true,
    grid style=dashed,
    legend style={fill=none},
    legend pos={south east},
    ytick = {1,3,5},
]
\addplot[color=black]
table[x index = 0,y index=1,col sep=comma]{quad_sum_wghts.tex};
\addplot[color=myRed]
table[x index = 0,y index=1,col sep=comma]{quad_sum_wghts_use_exact_int.tex};
\end{axis}
\end{tikzpicture}
\caption{$\sum_{\nu}|c_{\nu}|$ as a function of the requested tolerance for $\Psi_1$}
\label{subfig:wghts1}
\end{subfigure}
\hfill
\begin{subfigure}[t]{.49\linewidth}
    \begin{tikzpicture}
    \pgfplotsset{
  /pgfplots/xlabel near ticks/.style={
     /pgfplots/every axis x label/.style={
        at={(ticklabel cs:0.5)},anchor=near ticklabel
     }
  },
  /pgfplots/ylabel near ticks/.style={
     /pgfplots/every axis y label/.style={
        at={(ticklabel cs:0.5)},rotate=90,anchor=near ticklabel}
     }
  }
\begin{axis}[%
width=\linewidth,
    xmode=log,
    log basis x=10,
    xlabel={$\mathbf{tol}$},
    ylabel={$\sum_{\nu}|c_{\nu}|$},
    ylabel near ticks,
    ymajorgrids=true,
    xmajorgrids=true,
    grid style=dashed,
    legend style={fill=none},
    legend pos={south east}
]
\addplot[color=black]
table[x index = 0,y index=1,col sep=comma]{quad_sing_sum_wghts.tex};
\addplot[color=myRed]
table[x index = 0,y index=1,col sep=comma]{quad_sing_sum_wghts_use_exact_int.tex};
\end{axis}
\end{tikzpicture}
\caption{$\sum_{\nu}|c_{\nu}|$ as a function of the requested tolerance for $\Psi_2$}
\label{subfig:wghts2}
\end{subfigure}
    \caption{$\sum_{\nu}|c_{\nu}|$ for the QR-AAA approximation of $\Psi_1$ and $\Psi_2$, as a function of the requested tolerance. Version using exactly known integrals shown in red, Gauss-Kronrod based version shown in black.}
    \label{fig:wghts}
\end{figure}

While most are positive, some of the computed weights are negative. This is not out of the ordinary for general quadrature rules; in fact, as shown in~\cite{negWeights}, this can easily happen if the selected nodes are a large deviation of an equispaced grid. Additionally, the computed weights do have the desirable property that, for both strategies and both families of functions, the sum of the absolute values of the weights (see Figure~\ref{fig:wghts}) are well-bounded for $\mathbf{tol}\rightarrow 0$ and this bound improves as the number of quadrature nodes increases. This indicates that we have numerically stable quadrature rules. It is, however, still a limitation of the current version of our method; at this time we do not have a way to enforce positivity. From Figure~\ref{fig:poles} one can also see that the quadrature nodes cluster towards the endpoints of the domain, and they cluster stronger towards the singularity.

In conclusion, QR-AAA seems well-suited to construct quadrature rules for arbitrary classes of functions, provided that they can be well-approximated using rational functions with shared poles.

\subsection{Simple multivariate rational approximation}\label{sec:multivariate}
It is considered an important open problem to approximate multivariate rational functions by \textit{some} AAA-like algorithm. For simplicity we will restrict our exposition here largely to the bivariate case, but the multivariate case is completely similar.

Say we wish to approximate some function $f(x,y)$ by a rational function. Currently, two classes of methods considered to be state-of-the-art are the class of \textit{p-AAA type} methods (see \cite{pAAA,lrAAA}) and the class of \textit{parametric Loewner} methods (see \cite{pLoewner1,doi:10.1137/130914619}). In their basic form, similarly to univariate AAA, both construct approximants of $f(x,y)$ of the form
$$f(x,y)\approx r(x,y)=\left.\sum_{i=1}^{m_x}\sum_{j=1}^{m_y}\frac{w_{ij}f(x_i,y_j)}{(x-\xi_i)(y-\eta_j)}\middle/\sum_{i=1}^{m_x}\sum_{j=1}^{m_y}\frac{w_{ij}}{(x-\xi_i)(y-\eta_j)}\right. .$$
Parametric Loewner methods perform this approximation in a single step, by computing a null vector of a large Loewner matrix. The final approximation is then interpolatory in a predetermined set of points. 

By contrast, in p-AAA, the approximant is built up iteratively, by greedily selecting support points $(\xi_i,\eta_j)\in Z\subset \R^2$ from a set of samples $Z$ and determining the weights by minimizing a linearized least-squares error over the rest of the samples. This is done by computing a right singular vector associated with the minimal singular value of a large Loewner matrix. In its original form, the p-AAA algorithm was restricted to recovering functions $f(x,y)$ from tensorial sample grids $Z=Z_x\times Z_y$ of size $N_x \times N_y$, but in \cite{pAAAscattered} the method was extended to be applicable to scattered data. Even so, the tensor grid case, like in the case of multivariate polynomial approximation, is important enough to merit study.

For both p-AAA and parametric Loewner methods, the size of the Loewner matrices involved grows quickly as the number of dimensions $d$ grows. This is known as the \textit{curse of dimensionality}. For example, for ``vanilla'' p-AAA, applied to a $d$-variate function $f(x_1,\ldots,x_d)$, at each step a singular value decomposition of a Loewner matrix $\mathbb{L}\in \R^{\mathcal{O}(N^d)\times \mathcal{O}(m^d)}$ must be computed, where $N$ and $m$ are the number of sample points and the rational order per dimension respectively\footnote{For simplicity we have assumed an equal number of sample points and support points in each direction.}. Unlike in SV-AAA type methods, no simple SVD update scheme like in the SV-AAA case is currently available (see \cite{KarlSVAAA}). This means that at iteration $m$ the computational cost of p-AAA scales like $\mathcal{O}(m^{2d}N^{d})$. This quickly becomes intractable, especially in higher dimensions. 

A partial solution to the curse of dimensionality, low-rank AAA, was proposed in \cite{lrAAA}, where a low CPD tensor rank structure is imposed on the weight tensor, meaning $\mathbf{w} = \sum_{i=1}^{k}\mathbf{w}_{k,1}\otimes \cdots \otimes\mathbf{w}_{k,d}$. Another recent approach to tame the curse of dimensionality, called the \textit{Multivariate Loewner framework} was proposed in~\cite{MLFpaper}. There, an interpolatory approximation to a given multivariate function is built up using null spaces of $d$-dimensional Loewner matrices, which can be obtained from the null spaces of sequences of one-dimensional Loewner matrices.

We propose an alternate approach, based on the observation in Proposition~\ref{prop:2D}. Formulated informally:
\begin{proposition}\label{prop:2D}
    With the notation from Section~\ref{sec:qr-aaa}, suppose we apply QR-AAA on $f(x,y)$, sampled on $Z_x\times Z_y$, and obtain an interpolative decomposition 
    $$\mathbf{F}(Z_x,Z_y) \approx \mathbf{H}_x \mathbf{F}(\{\xi_1,\ldots,\xi_{m_x}\},Z_y)$$
    with weights $\{w_{i}\}_{i=1}^{m_x}$. Suppose by applying QR-AAA on $f(y,x)$ we similarly obtain an interpolative decomposition 
    $$\mathbf{F}(Z_x,Z_y) \approx \mathbf{F}(Z_x,\{\eta_1,\ldots,\eta_{m_y}\})\mathbf{H}_y^T$$
    with weights $\{w_{j}\}_{j=1}^{m_y}$.
    Then
    $$\mathbf{F}(Z_{x},Z_y) \approx \mathbf{H}_x \mathbf{F}(\{\xi_1,\ldots,\xi_{m_x}\},\{\eta_1,\ldots,\eta_{m_x}\})\mathbf{H}_y^T$$
    and this expression can be written explicitly as
    \begin{equation}\label{eq:repSepAAA}
        f(x,y)\approx r(x,y)=\left.\sum_{i=1}^{m_x}\sum_{j=1}^{m_y}\frac{w_{i}w_{j}f(x_i,y_j)}{(x-\xi_i)(y-\eta_j)}\middle/\sum_{i=1}^{m_x}\sum_{j=1}^{m_y}\frac{w_{i}w_{j}}{(x-\xi_i)(y-\eta_j)}\right.
    \end{equation}
    for $(x,y)\in Z_x\times Z_y$.
\end{proposition}
This will not be proven here, but it should be mentioned that it is essential (at least for the theory) that the set-valued approximations are applied to stable bases that span the rows and columns of $\mathbf{F}(Z_x,Z_y)$, and that the matrices $\mathbf{H}_x,\mathbf{H}_y$ are ``well-behaved'' in the sense of Section~\ref{sec:qr-aaa}.

Note the similarity to the low-rank p-AAA variant. In fact, the \textit{representation} in equation~\eqref{eq:repSepAAA} corresponds to the rank one case. The main difference between the method from~\cite{lrAAA} and our method lies in how the support points and weights are obtained. Here, the variables are first treated completely independently, after which the QR-AAA approximations are joined. As will be shown in Section~\ref{subsubsec:numExpMulti} this leads to a fast and effective scheme for a large class of multivariate functions. 

In the multivariate case, the two-sided interpolative decomposition in equation~\eqref{eq:repSepAAA} becomes an interpolative Tucker decomposition:
\begin{equation}\label{eq:tmprod}
    \mathbf{F}(Z_1,\ldots,Z_d)\approx F(\Xi_1,\ldots,\Xi_d)\bigtimes_{\ell=1}^d\mathbf{H}_{\ell}
\end{equation}
in which $\Xi_{\ell}$ and $\mathbf{H}_{\ell}$ are the mode $\ell$ support nodes and barycentric matrix respectively. The symbol `$\bigtimes_{\ell=1}^d$' denotes the sequence of mode $\ell$ tensor-matrix products, as $\ell$ ranges from $1$ to $d$.

Our method of multivariate approximation summarized in Algorithm~\ref{alg:Multivariate}. We use \texttt{MATLAB} notation wherever convenient.

\begin{algorithm}[h]
\caption{Multivariate rational approximation by repeated QR-AAA}\label{alg:Multivariate}
\SetKwInOut{Input}{Input}
\SetKwInOut{Output}{Output}
\SetKw{Init}{init}{}{}
\SetAlgoLined
\DontPrintSemicolon
\Input{Tensor $\mathbf{F}\in\mathbb{R}^{N_1\times\cdots\times N_d}$ and associated sample points $(Z_1\times\cdots\times Z_d)$, tolerances \texttt{tol\_qr}, \texttt{tolAAA}.}
\Output{Multivariate rational approximant $r(x_1,\ldots,x_d)$.}
 \For{$\ell=1\ldots d$}{
    \tcp{$\mathbf{F}_{\ell}$ is mode $\ell$ unfolding:}
    $[\mathbf{Q},\mathbf{R},\sim] = \texttt{qr}(\mathbf{F}_{\ell},\texttt{`econ'})$\;
    \texttt{dr} $=$ \texttt{abs(diag(}$\mathbf{R}$\texttt{))}\;
    \texttt{rk} $=$ \texttt{sum(dr>tol\_qr*dr(1))}\;
    $\mathbf{Q}_{\ell} = \mathbf{Q}$\texttt{(:,rk).*dr(1:rk).'}\;
    \tcp{Vector-valued rational function, nodes, weights: }
    \texttt{[$r_{\ell}$, $\Xi_{\ell}$, $\mathbf{w}_{\ell}$]} $=$ \texttt{aaa\_sv($\mathbf{Q}_{\ell}$,$Z_{\ell}$,`tol',tolAAA)}\;
 }
 \tcp{With `r' corresponding to equation~\eqref{eq:repSepAAA}:}
    $r := $\texttt{@($x_1,\ldots,x_d$)} \texttt{r($\{\Xi_{\ell}\}_{\ell}$,$\{\mathbf{w}_{\ell}\}_{\ell}$,$F(\Xi_1,\ldots,\Xi_d)$)($x_1,\ldots,x_d$)}\;
\end{algorithm}

Alternatively, one can use the component-wise QR decompositions\footnote{Or really any method to compute a stable tensor ID $\mathbf{F}\approx \mathbf{F}(J_1,\ldots,J_{d})\bigtimes_{\ell=1}^d\mathbf{U}_{\ell}$.} only to construct the indices $\{J_{\ell}\}_{\ell=1}^d$ in an interpolative decomposition for $\mathbf{F}$, after which SV-AAA is applied to the component functions $\mathbf{F}(J_1,J_2,\ldots,J_{\ell-1},:,J_{\ell+1},\ldots,J_{d})$ for $\ell=1,\ldots,d$. Both approaches are included in our \texttt{MATLAB} code, but since the results do not differ significantly, we only report on the first.

The data-efficiency of the final rational representation depends only on the degrees of the component-wise vector-valued rational approximations. For sufficiently smooth multivariate functions, it is to be expected that these degrees satisfy $m_{\ell}\ll N_{\ell}$, and good accuracy can be reached even by low-degree rational approximants. 

We will see in Section~\ref{subsubsec:twoStep} that adversarial examples do exist. This happens if there is no ``good'' vector-valued rational approximant of (one of) the component functions $x_i\mapsto f(Z_1,\ldots,x_i,\ldots,Z_d)$. As a partial remedy, a \textit{two-step procedure} is introduced in Section~\ref{subsubsec:twoStep}

\subsubsection{Numerical experiments}\label{subsubsec:numExpMulti}
In this section, the QR-AAA based multivariate rational approximation scheme outlined in Algorithm~\ref{alg:Multivariate} is applied to $4$ functions, representing some important functions classes, all sampled on uniform grids. These are:
\begin{enumerate}
    \item $f_1(x,y,z) = \frac{1}{xyz+2}$ with $x,y,z\in[-1,1]$: a simple, non-separable, rational function. Here $N_x=N_y=N_z=51$.
    \item $f_2(x,y,z) = \frac{1}{\sqrt{x^2+2y^2+3z^2+\delta}}$ for $\delta\rightarrow 0^+$ with $x,y,z\in [-1,1]$: a non-rational function with a singularity approaching the domain. Here $N_x=N_y=N_z= 151$.
    \item $f_3(x,y,z) = \sqrt{x^2+2y+3z}$ with $x\in [-1,1]$ and $y,z\in[0,1]$: a function containing multiple types of singularity; algebraic in $y,z=0$ and $C^0$-type in $x=0$\footnote{By this we mean that around $0$, in $x$, the function is $C^0$, but not $C^1$}. Here $N_x=401$ $N_y=N_z= 200$.
    \item $f_4(x,y,z,t) = \frac{\cos(3\pi x)\sin(4\pi y)+\sin(2\pi z)\sin(5\pi t)}{xy^2+2z^3t+4}$ for $x,y,z,t\in [-1,1]$: a 4-D oscillatory function multiplied by a degree $4$ rational function. Here $N_x=N_y=N_z=N_t=51$.
\end{enumerate}
In examples $2$ through $4$, the coefficients for the variables are deliberately chosen to be different (as opposed to choosing, say, $f(x,y,z)=1/\sqrt{x^2+y^2+z^2+\delta}$). This was done to break symmetry as much as possible.
\begin{figure}
    \centering
    \begin{subfigure}{.46\linewidth}
        \begin{tikzpicture}
            \begin{axis}[%
                width=\linewidth,
                ymode=log,
                log basis y=10,
                xlabel={degree $m$},
                ylabel={$\|\mathbf{f}-\mathbf{r}\|_F/\|\mathbf{f}\|_F$},
                ylabel near ticks,
                ymajorgrids=true,
                xmajorgrids=true,
                grid style=dashed
                ]
                \addplot[color=black]
                table[x index = 0,y index=1,col sep=comma]              {conv_mmax_grid_ratfunc.tex};
                \addlegendentry{$Z$}
                \addplot[color=black,dashed]
                table[x index = 0,y index=1,col sep=comma]{conv_mmax_off_grid_ratfunc.tex};
                \addlegendentry{$Z_{\text{val}}$}
            \end{axis}
        \end{tikzpicture}
        \caption{Relative Frobenius norm error for $f_1(x,y,z)$ on sample grid $Z$ and validation grid $Z_{\text{val}}$.}
        \label{subfig:convmmax1}
    \end{subfigure}
    \hfill
    \begin{subfigure}{.46\linewidth}
        \begin{tikzpicture}
            \begin{axis}[%
                width=\linewidth,
                ymode=log,
                log basis y=10,
                xlabel={degree $m$},
                ylabel={$\|\mathbf{f}-\mathbf{r}\|_F/\|\mathbf{f}\|_F$},
                ylabel near ticks,
                ymajorgrids=true,
                xmajorgrids=true,
                grid style=dashed
                ]
                \addplot[color=black]
                table[x index = 0,y index=1,col sep=comma]              {conv_mmax_grid_invsqrt.tex};
                \addlegendentry{$Z$}
                \addplot[color=black,dashed]
                table[x index = 0,y index=1,col sep=comma]{conv_mmax_off_grid_invsqrt.tex};
                \addlegendentry{$Z_{\text{val}}$}
            \end{axis}
        \end{tikzpicture}
        \caption{Relative Frobenius norm error for $f_2(x,y,z)$ on sample grid $Z$ and validation grid $Z_{\text{val}}$.}
        \label{subfig:convmmax2}
    \end{subfigure}
    \\
    \begin{subfigure}{.46\linewidth}
        \begin{tikzpicture}
            \begin{axis}[%
                width=\linewidth,
                ymode=log,
                log basis y=10,
                xlabel={degree $m$},
                ylabel={$\|\mathbf{f}-\mathbf{r}\|_F/\|\mathbf{f}\|_F$},
                ylabel near ticks,
                ymajorgrids=true,
                xmajorgrids=true,
                grid style=dashed
                ]
                \addplot[color=black]
                table[x index = 0,y index=1,col sep=comma]              {conv_mmax_grid_sqrt.tex};
                \addlegendentry{$Z$}
                \addplot[color=black,dashed]
                table[x index = 0,y index=1,col sep=comma]{conv_mmax_off_grid_sqrt.tex};
                \addlegendentry{$Z_{\text{val}}$}
            \end{axis}
        \end{tikzpicture}
        \caption{Relative Frobenius norm error for $f_3(x,y,z)$ on sample grid $Z$ and validation grid $Z_{\text{val}}$.}
        \label{subfig:convmmax3}
    \end{subfigure}
    \hfill
    \begin{subfigure}{.46\linewidth}
        \begin{tikzpicture}
            \begin{axis}[%
                width=\linewidth,
                ymode=log,
                log basis y=10,
                xlabel={degree $m$},
                ylabel={$\|\mathbf{f}-\mathbf{r}\|_F/\|\mathbf{f}\|_F$},
                ylabel near ticks,
                ymajorgrids=true,
                xmajorgrids=true,
                grid style=dashed
                ]
                \addplot[color=black]
                table[x index = 0,y index=1,col sep=comma]              {conv_mmax_grid_oscil.tex};
                \addlegendentry{$Z$}
                \addplot[color=black,dashed]
                table[x index = 0,y index=1,col sep=comma]{conv_mmax_off_grid_oscil.tex};
                \addlegendentry{$Z_{\text{val}}$}
            \end{axis}
        \end{tikzpicture}
        \caption{Relative Frobenius norm error for $f_4(x,y,z,t)$ on sample grid $Z$ and validation grid $Z_{\text{val}}$.}
        \label{subfig:convmmax4}
    \end{subfigure}
    \caption{Convergence plot, over the maximal degree $m:=\max_{\ell}\{m_\ell\}$, for the $4$ functions considered. For $f_2$, $\delta=2^{-6}$ is used.}
    \label{fig:QRAAAmultiResults}
\end{figure}
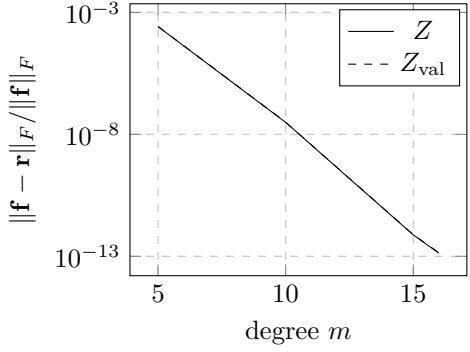
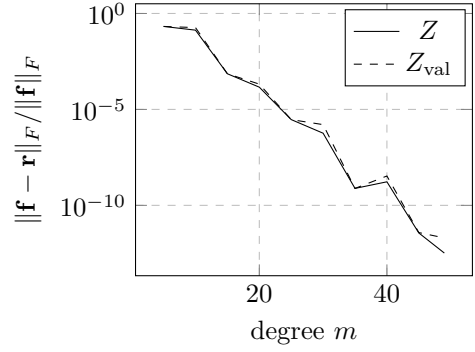
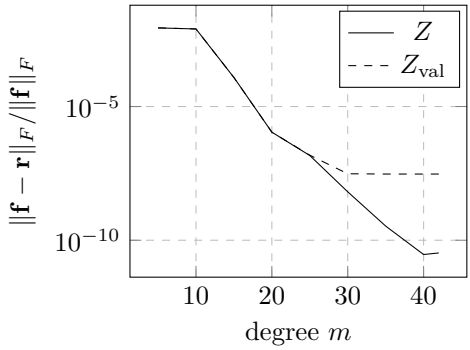
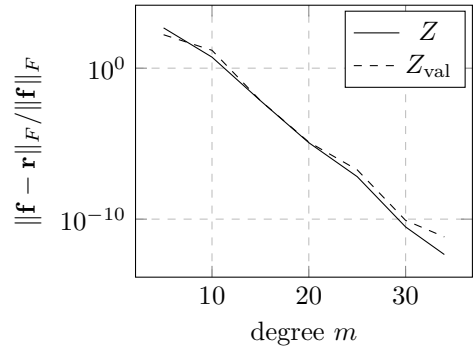
The results are reported in Fgure~\ref{fig:QRAAAmultiResults}. To be explicit, the validation grid in each case was chosen to be a finer grid on the same domain, and the error studied is the Frobenius norm error over this finer tensor grid. 

We see that for all $4$ functions, excellent convergence on the sample grid is obtained. However, for $f_3$, the validation error stalls around $2\cdot10^{-8}$. This hints at the fact that the sample grid is not well-chosen; as in univariate rational approximation to the square root, uniform support points are not optimal. Which multidimensional grid would be optimal is an open question.

\begin{figure}
    \centering
    \begin{subfigure}{.46\linewidth}
        \begin{tikzpicture}
            \begin{axis}[%
                width=\linewidth,
                xmode=log,
                log basis x=2,
                ymode=log,
                log basis y=10,
                xlabel={$\delta$},
                ylabel={$\|\mathbf{f}-\mathbf{r}\|_F/\|\mathbf{f}\|_F$},
                ylabel near ticks,
                ymajorgrids=true,
                xmajorgrids=true,
                grid style=dashed
                ]
                \addplot[color=black]
                table[x index = 0,y index=1,col sep=comma]{err_del_grid_invsqrt.tex};
                \addlegendentry{$Z$}
                \addplot[color=black,dashed]
                table[x index = 0,y index=1,col sep=comma]{err_del_off_grid_invsqrt.tex};
                \addlegendentry{$Z_{\text{val}}$}
            \end{axis}
        \end{tikzpicture}
        \caption{Error for $f_2(x,y,z)$ on sample grid $Z$ and validation grid $Z_{\text{val}}$, as $\delta\rightarrow 0^+$.}
        \label{subfig:delerr1}
    \end{subfigure}
    \hfill
    \begin{subfigure}{.46\linewidth}
        \begin{tikzpicture}
            \begin{axis}[%
                width=\linewidth,
                xlabel={$\delta$},
                ylabel={$m$},
                xmode = log,
                log basis x=2,
                ylabel near ticks,
                ymajorgrids=true,
                xmajorgrids=true,
                grid style=dashed
                ]
                \addplot[color=black]
                table[x index = 0,y index=1,col sep=comma]{degr_del_invsqrt.tex};
                \addlegendentry{$m_x$}
                \addplot[color=black,dashed]
                table[x index = 0,y index=2,col sep=comma]{degr_del_invsqrt.tex};
                \addlegendentry{$m_y$}
                \addplot[color=black,dotted]
                table[x index = 0,y index=3,col sep=comma]{degr_del_invsqrt.tex};
                \addlegendentry{$m_z$}
            \end{axis}
        \end{tikzpicture}
        \caption{Degrees $m_x,m_y,m_z$ for the rational approximant to $f_2(x,y,z)$ as $\delta\rightarrow 0^+$}
        \label{subfig:degrdel}
    \end{subfigure}
    \caption{Relative Frobenius norm error and rational degrees, as $\delta\rightarrow 0^+$, for $f_2(x,y,z)$.}
    \label{fig:QRAAAdelResults}
\end{figure}
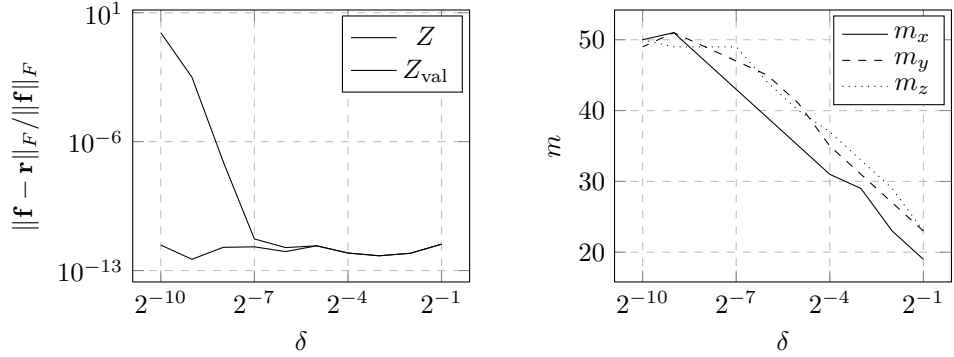

For $f_2(x,y,z)$, a study for $\delta\rightarrow 0^+$ is provided in Figure~\ref{fig:QRAAAdelResults}. As the singularity approaches the domain, the rational degrees increase linearly, up to the point where $m_x=N_x/2$, $m_y=N_y/2$ and $m_z=N_z/2$, at which point all accuracy is lost on the validation grid. Indeed, at this point, the degree needed to accurately resolve the near-singularity cannot be achieved with the provided number of sample points.

To close this section, in Table~\ref{tab:statsMulti}, the timings and degrees at a fixed tolerance of $10^{-8}$ for $f_1-f_4$ are included. Let us stress that this is using \textit{full} QR-decompositions,  as \texttt{MATLAB} currently does not support an early-exit Householder column-pivoted QR decomposition. Using the implementation (serial or parallel) from \cite{PQRAAA} would drastically improve the timings for the QR-step. 

\begin{table}[]
    \centering
    \begin{tabular}{c|c|c|c|c|c}
    &$\mathbf{err}_{\text{val}}$&$t_{\text{QR}}$&$t_{\text{AAA}}$&ranks&degrees\\
    \hline
    $f_1(x,y,z)$&$4.7\cdot10^{-9}$&$0.01$s&$6$ms&$(12,12,12)$&$(10,10,10)$\\
    $f_2(x,y,z)$&$3.13\cdot 10^{-9}$&$2.2$s&$30$ms&$(15,16,16)$&$(30,32,32)$\\
    $f_3(x,y,z)$&$1.6\cdot 10^{-8}$&$19.1$s&$35$ms&$(13,13,13)$&$(33,13,13)$\\
    $f_4(x,y,z,t)$&$2.39\cdot 10^{-8}$&$3.63$s&$88$ms&$(14,14,16,18)$&$(16,23,19,21)$
    \end{tabular}
    \caption{Relative Frobenius norm errors on a validation grid, timings for the QR-decompositions and the vector-valued AAA approximations, QR ranks and AAA degrees at a fixed tolerance $10^{-8}$ for the multivariate functions $f_1-f_4$. For $f_2(x,y,z)$, $\delta$ was set to $2^{-6}$.}
    \label{tab:statsMulti}
\end{table}

Observe that the rational degree for the approximation of $f_1$, in this case $10$, exceeds the actual rational degree, which is $1$. This is in contrast to p-AAA type methods and parametrized and multivariate Loewner methods, which recover the exact rational degree of a given rational function. Additionally, when supplied with a rational function, these methods recover the function to machine precision. For irrational functions, however, the picture is different. In Figure~\ref{fig:compareMLFpAAA}, we compare as a function of the approximation degree, the error and computation time for
\begin{enumerate}
    \item our proposed approach, and
    \item the low-rank AAA method.
\end{enumerate}
We consider only the error on a validation grid. Due to the high computation time for low-rank AAA, we have reduced the grid to $N_x = N_y = N_z=75$.
\begin{remark}
    For the irrational functions considered here, the multivariate Loewner Framework (mLf) methods did not produce convergent approximations. The direct version in ~\cite[Algorithm 7.1]{MLFpaper} does attain $\sim 2$ digits of accuracy, but does not detect higher degree is needed for a more accurate approximation. The adaptive version in~\cite[Algorithm 7.2]{MLFpaper} does not converge to the underlying function on the validation grid. However, there are many irrational functions for which the mLf methods \textbf{do} perform well, as shown in~\cite{mlfCompare}.
\end{remark}
\begin{figure}
    \centering
    \begin{subfigure}{.46\linewidth}
        \begin{tikzpicture}
            \begin{axis}[%
                width=\linewidth,
                ymode=log,
                log basis y=10,
                xlabel={$\max.\textrm{degr.}$},
                ylabel={$\|\mathbf{f}-\mathbf{r}\|_F/\|\mathbf{f}\|_F$},
                ylabel near ticks,
                ymajorgrids=true,
                xmajorgrids=true,
                grid style=dashed,
                legend style={font=\footnotesize}
                ]
                \addplot[color=black,mark=*]
                table[x index = 0,y index=1,col sep=comma]{err_off_grid.tex};
                \addlegendentry{QR}
                \addplot[color=black]
                table[x index = 0,y index=1,col sep=comma]{err_off_grid_paaa.tex};
                \addlegendentry{LR}
            \end{axis}
        \end{tikzpicture}
        \caption{Error for $f_2(x,y,z)$, $\delta=2^{-6}$, on the validation grid $Z_{\rm val}$.}
        \label{subfig:delerr2}
    \end{subfigure}
    \hfill
    \begin{subfigure}{.46\linewidth}
        \begin{tikzpicture}
            \begin{axis}[%
                width=\linewidth,
                ymode=log,
                log basis y=10,
                xlabel={$\max.\textrm{degr.}$},
                ylabel={Time [s]},
                ylabel near ticks,
                ymajorgrids=true,
                xmajorgrids=true,
                grid style=dashed,
                legend style={font=\footnotesize},
                legend pos ={north west}
                ]
                \addplot[color=black,mark=*]
                table[x index = 0,y index=1,col sep=comma]{time_qraaa.tex};
                \addlegendentry{QR}
                \addplot[color=black]
                table[x index = 0,y index=1,col sep=comma]{time_paaa.tex};
                \addlegendentry{LR}
            \end{axis}
        \end{tikzpicture}
        \caption{Timings for $f_2(x,y,z)$, $\delta=2^{-6}$.}
        \label{subfig:timeComparepAAA}
    \end{subfigure}
    \caption{Relative Frobenius norm error and timing over the rational degrees for our QR-AAA method (QR) and the low-rank AAA method (LR).}
    \label{fig:compareMLFpAAA}
\end{figure}

\subsubsection{A two-step procedure}\label{subsubsec:twoStep}
There are at least two scenarios in which we would want to augment the method from Algorithm~\ref{alg:Multivariate}:

\begin{enumerate}
    \item There is at least one direction for which QR-AAA fails. In practice this can be detected when the Loewner matrices in the QR-AAA algorithm become singular, or the desired precision is not reached before the maximal degree is reached
    \item We explicitly \textit{do not want} a separable approximation to $f(x_1,\ldots,x_n)$. This can be particularly useful in the case of \textit{function extension}, as we will see in Section~\ref{sec:extension}.
\end{enumerate}

We propose the following simple modification:
\begin{enumerate}
    \item First, run Algorithm~\ref{alg:Multivariate} with a fixed (low) tolerance or maximum degree for the desired direction(s).
    \item Subsequently, use the selected support points as initial support points for p-AAA.
\end{enumerate}
By this we mean that we force the interpolation points in the p-AAA algorithm to include the selected points, after which the usual p-AAA procedure starts.

As an illustrating example, we consider a synthetic rational transfer function $H(s,p)$ from~\cite{lrAAA}, with $50$ poles close to the domain of approximation. It is sampled on $s\in[\imath,10^{4}\imath ]$ and $p\in[10^{-1.5},1]$, both being discretized using log-spaced grids ($N_s=500,N_p=50$). For this function, it is only the $s$-direction for which the QR-AAA approach fails. We therefore set the maximal QR-AAA degree in the $s$-direction to be $30$, and simultaneously set the desired accuracy for the QR-AAA approximation to be $10^{-5}$ in the $p$-direction\footnote{This is to reduce $m_p$ at the start of the p-AAA step. While the two-step method still converges in case the degree is overestimated, recall that the sizes of the Loewner matrices for p-AAA grow like $\mathcal{O}(N_sN_p)\times\mathcal{O}(m_sm_p)$}. 
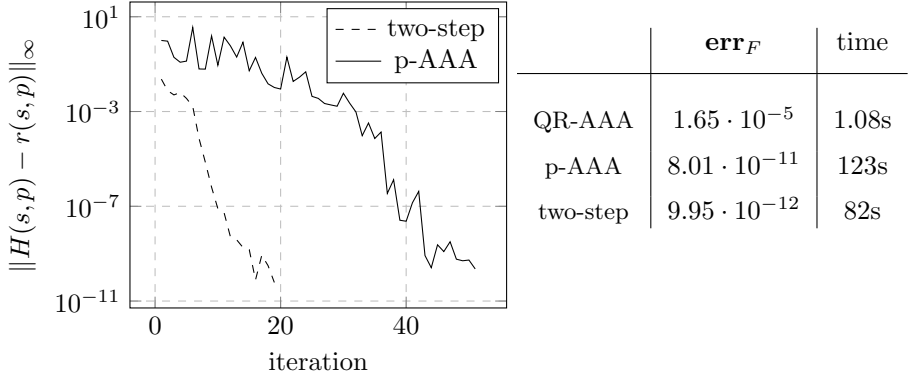
\begin{figure}
    \centering
    \begin{subfigure}[b]{.5\linewidth}
    \begin{tikzpicture}
    \begin{axis}[
    width=\linewidth,
    ymode=log,
    log basis y=10,
    xlabel={iteration},
    ylabel={$\|H(s,p)-r(s,p)\|_{\infty}$},
    ymajorgrids=true,
    xmajorgrids=true,
    grid style=dashed
]
\addplot[color=black,dashed] 
    table [x index = 0,y index = 1,col sep=comma]{conv_2step.tex};
    \addlegendentry{two-step}
\addplot[
    color=black
    ] table [x index = 0,y index = 1,col sep=comma]{conv_pAAA.tex};
    \addlegendentry{p-AAA}
    \end{axis}
    \end{tikzpicture}
    \caption{Convergence in Frobenius norm error over refined validation grid}
    \label{subfig:HspConv}
    \end{subfigure}
    \hfill
    \begin{subfigure}[b]{.48\linewidth}
        \begin{tabular}{c|c|c}
    &$\mathbf{err}_{F}$&time\\[5pt]
    \hline
    &&\\
    \small QR-AAA&$1.65\cdot 10^{-5}$&$1.08$s\\[5pt]
    \small p-AAA&$8.01\cdot 10^{-11}$&$123$s\\[5pt]
    \small two-step&$9.95\cdot 10^{-12}$&$82$s
    \end{tabular}
    \vspace{2cm}
    \caption{Timings and Frobenius norm error in sample grid of the final approximation}
    \label{tab:Hspresults}
    \end{subfigure}
    
    \caption{Comparison of QR-AAA, p-AAA and two-step AAA for $H(s,p)$ the degree $50$ synthetic rational transfer function. Convergence plot~\ref{subfig:HspConv} shows the relative Frobenius norm error over a refined validation grid as a function of the p-AAA iteration number. For two-step AAA this is after the initial QR-AAA based component-wise approximations have been constructed.}
    \label{fig:NumericalResultsHsp}
\end{figure}

The results are reported in Figure~\ref{fig:NumericalResultsHsp}. The final degrees of the QR-AAA, p-AAA and two-step approximation were, respectively, $(100,32)$, $(49,19)$ and $(48,18)$.

\subsection{Multivariate function extension}\label{sec:extension}
Consider a real-analytic function $f(x_1,\ldots,x_n)$, given in some simply connected domain $\Omega\subset\mathbb{R}^n$. We consider here the problem of \textit{function extension}: based on samples of $f$ in $\Omega$, the aim of function extension is to construct an analytic function $\tilde{f}$ in $\widetilde{\Omega}\supset\Omega$ that ``agrees'' with $f$. Depending on the application, this may be in the sense that $\tilde{f}$ agrees with $f$ in $\widetilde{\Omega}$ up to some precision, or only that $\|f-\tilde{f}\|_{\infty}$ is sufficiently small in $\Omega$ and $\tilde{f}$ is analytic in $\widetilde{\Omega}$. We restrict our attention here to $n=2$, the bivariate case. As in Section~\ref{sec:multivariate}, however, the ideas outlined in this section generalize easily to higher dimensions.

\begin{figure}
    \centering
    \begin{subfigure}[t]{.32\linewidth}
    \includegraphics[width=\linewidth]{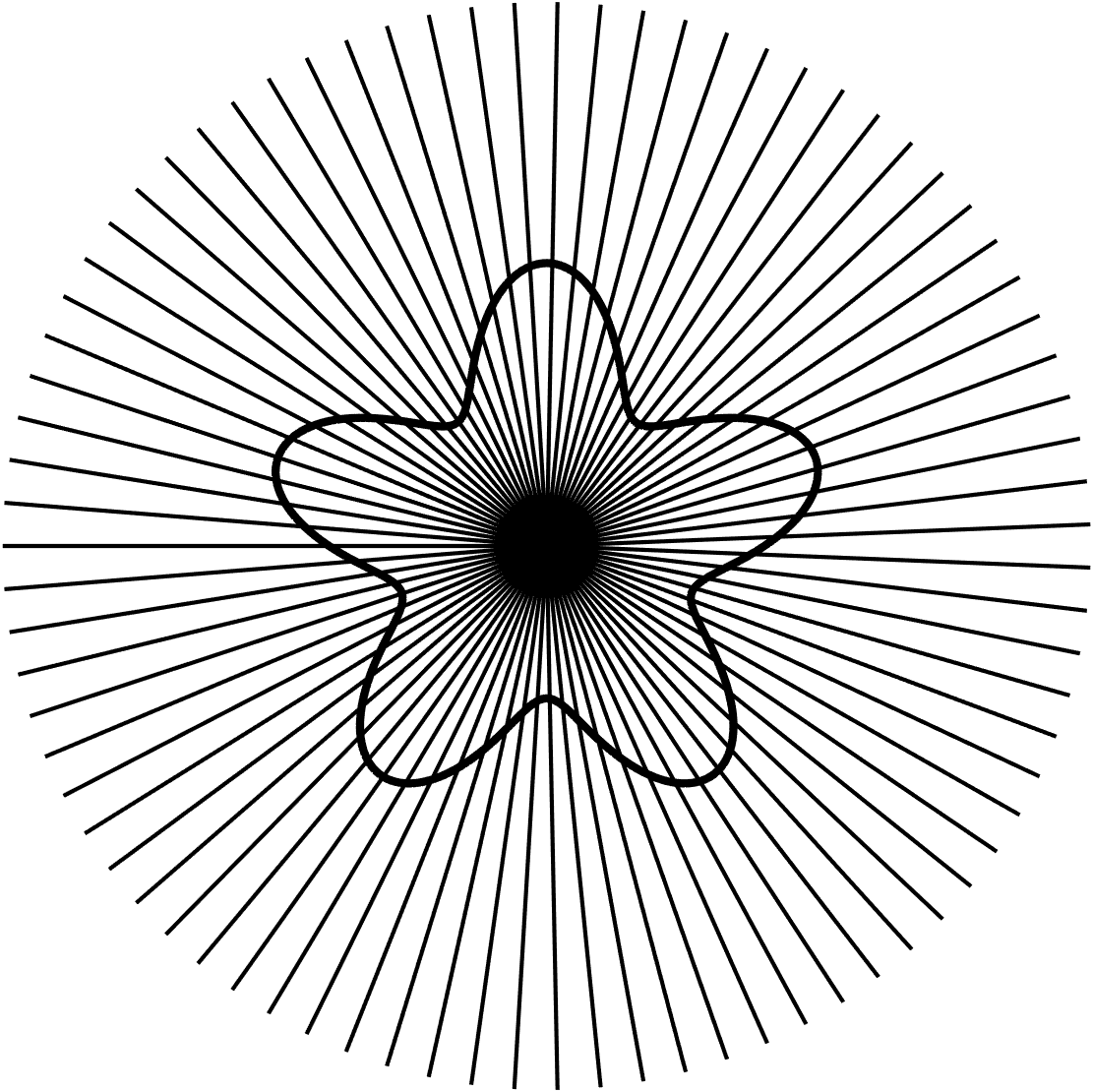}
        \caption{Set-up of the rays in ray-AAA for $\Omega$ a starfish domain.}
        \label{subfig:rays1}
    \end{subfigure}
    \hfill
    \begin{subfigure}[t]{.32\linewidth}
    \includegraphics[width=\linewidth]{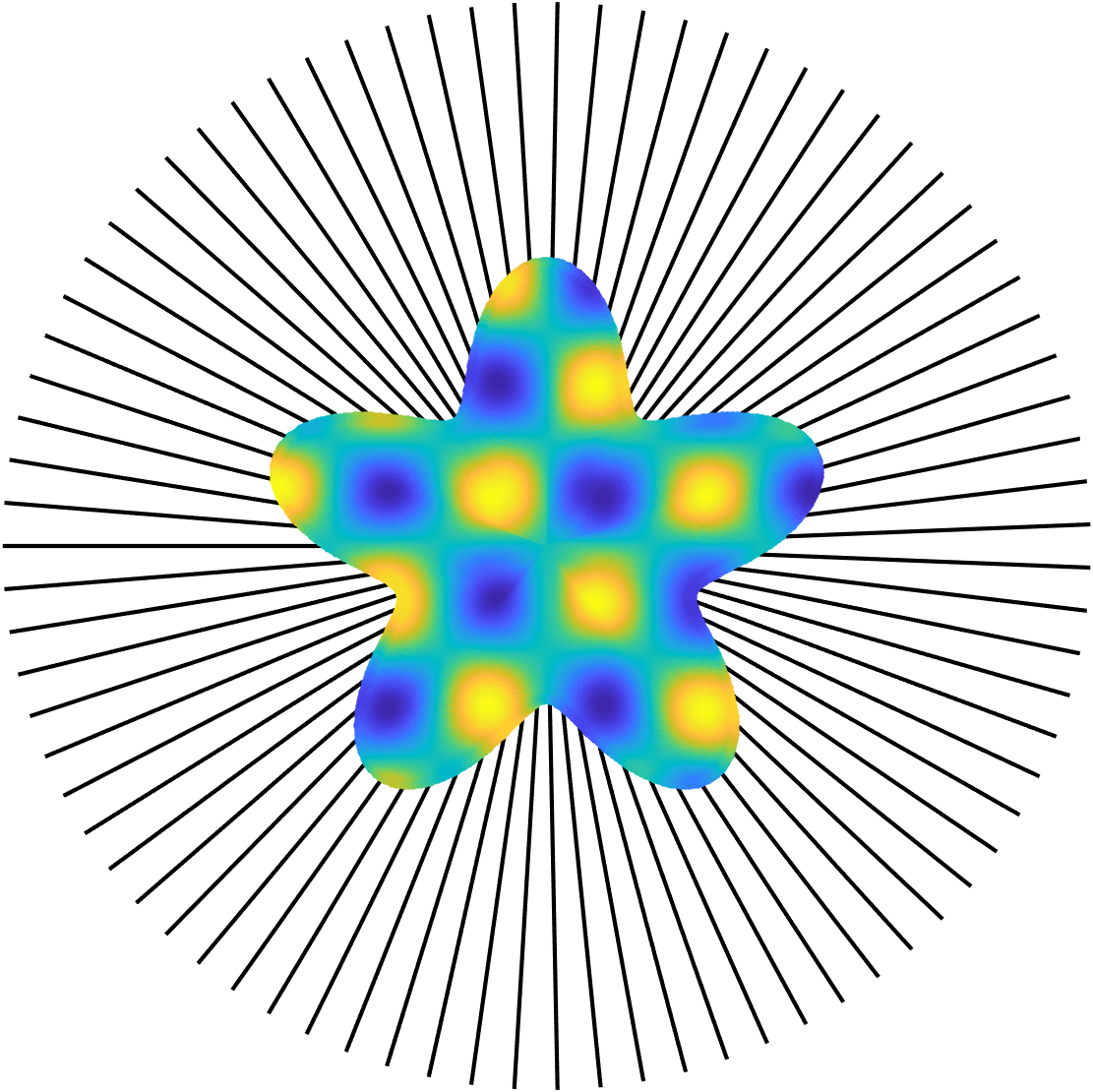}
        \caption{Function $f(x,y)$ to be extended, sampled in $\Omega$}
        \label{subfig:rays2}
    \end{subfigure}
    \hfill
    \begin{subfigure}[t]{.32\linewidth}
    \includegraphics[width=\linewidth]{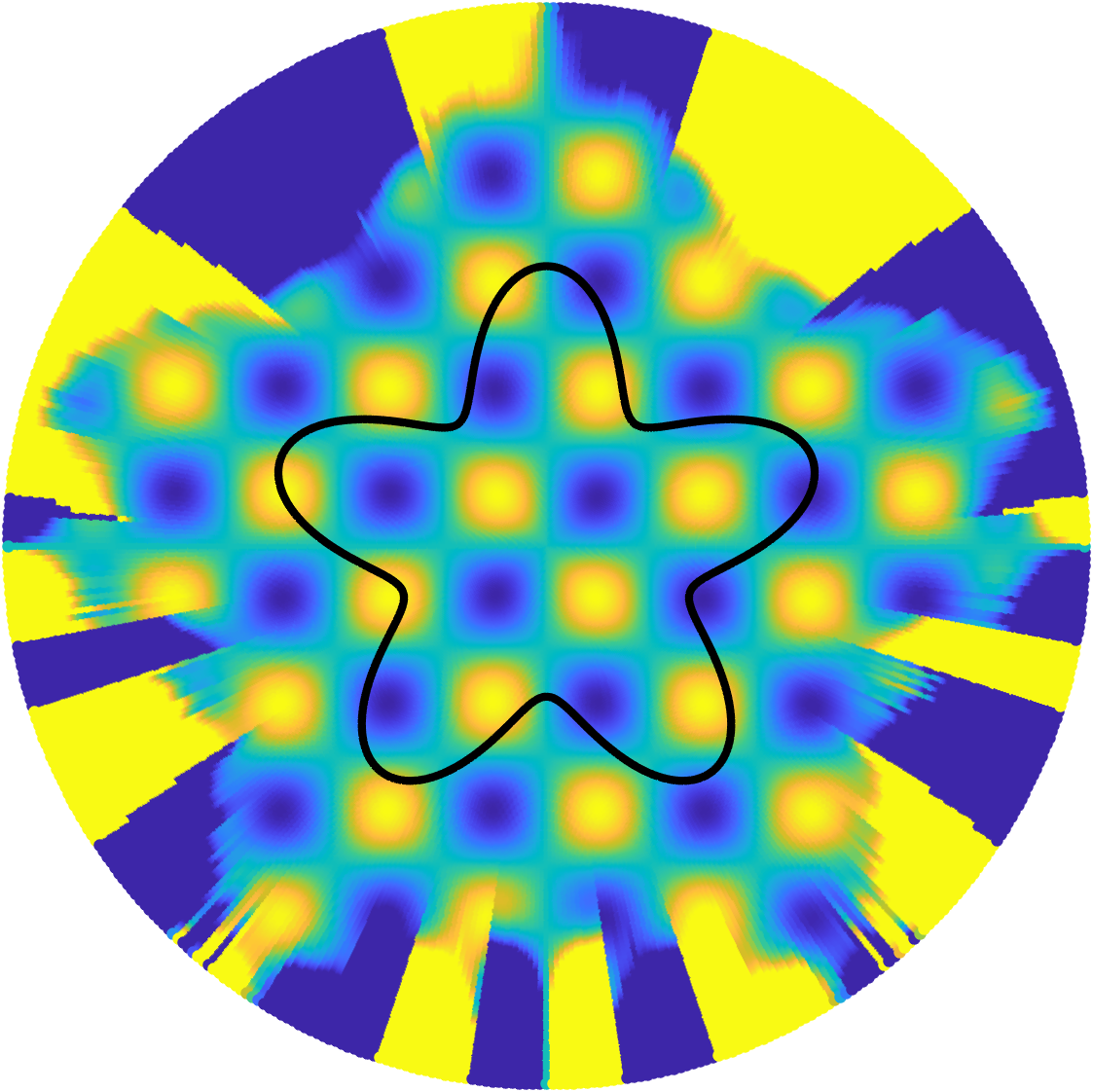}
        \caption{Extension obtained by the ray-AAA approximation}
        \label{subfig:extension}
    \end{subfigure}
    \caption{Principle of ray-AAA and plot of the ray-AAA extension of $f(x,y)=-\sin(2\pi x)\sin(2\pi y)$. In Figure~\ref{subfig:extension}, for convenience of plotting, the extension is set to $\pm 1$ when it becomes too large to visualize.}
    \label{fig:rayAAA}
\end{figure}

In \cite{NickExtension}, it was demonstrated that AAA rational approximation is an excellent tool for analytic continuation and function extension. In a nutshell, if a function $f(x,y)$ is to be extended from $\Omega\subset\mathbb{R}^2$ along some direction, $f(x,y)$ is sampled along a ray in this direction contained in $\Omega$. Subsequently, a rational approximation $r(x,y)$ to $f(x,y)$ is constructed from these sample points by means of AAA, and it is this rational function that is extended outside of $\Omega$. This process is repeated for each direction of interest. For the purposes of this section, this process will be referred to as \textit{ray-AAA}. An illustration of the ray-AAA method is shown in Figure~\ref{fig:rayAAA} (recreated from \cite{NickExtension}).

Immediately the similarity to the problem in Section~\ref{sec:Hedgehog} should be evident. Though there, through the addition of the on-surface points, the problem is changed from an extension problem to an interpolation problem, it is not unreasonable to assume that QR-AAA can again be used to improve the performance of ray-AAA.

Two problems will be considered:
\begin{enumerate}
    \item Extending the function $f(x,y) = -\sin(2\pi x)\sin(2\pi y)$ given on a starfish domain $\Omega$, whose boundary is parametrized by $\partial\Omega=\{\gamma(\theta):\theta\in[-\pi,\pi)\}$, where
    $$\gamma(\theta)=(1+0.3\sin(5\theta))\cdot(\cos(\theta),\sin(\theta)),$$
    i.e., the domain shown in Figure~\ref{fig:rayAAA}. This forms a good test case for \textit{oscillatory function extension}, which is known to be hard. In particular, as in \cite{NickExtension}, rational approximants in 16 digit arithmetic tend to only retain any accuracy up to one wavelength of extension. This is what is known as the \textit{one-wavelength principle}. Beyond this point, at least typically, all accuracy is lost.
    \item Extending the norm function $f(x,y)=\sqrt{x^2+y^2}$, given on the domain $\Omega = [-4,-2]\times[-4,4]$. Of particular interest here is the extension along the $x$-direction, where for $(x,y)=(0,0)$ the function $f(x,y)$ contains an algebraic singularity. Again, this is a challenging problem; in general it cannot be expected that an analytic extension, by its very nature, can model such a singularity.
\end{enumerate}
For each of these, a specialized method utilizing QR-AAA will be constructed. These are summarized in Section~\ref{sec:starfish} and Section~\ref{sec:root} respectively. 
\subsubsection{Oscillatory function extension}\label{sec:starfish}
Since the starfish domain $\Omega$ is periodic, one can write
$$f(x,y) = f(\gamma(\rho,\theta)):=f(\rho,\theta)$$
in which $\gamma(\rho,\theta) = \gamma(\theta)\rho$. The Greek letter $\rho$ is used here for the radial variable in place of $r$, since $r$ is reserved for rational approximants. For each $\rho_0\in[0,1]$ we make the ansatz that the function $f(\rho_0,\theta)$ is approximately a bandlimited function:
$$f(\rho_0,\theta)\approx\sum_{k=-m}^{m}c_k(\rho_0) e^{\imath k \theta}$$
with $m$ some predetermined maximal frequency. This approximation problem is also called trigonometric interpolation, and can easily be achieved using the DFT if equispaced rays $\{\theta_i\}_{i=1}^N$ are chosen. Even in the non-equispaced case, so long as $m$ is not too large, this step poses no significant challenge.

\begin{figure}
    \centering
    \begin{subfigure}[t]{.32\linewidth}
    \includegraphics[width=.85\linewidth]{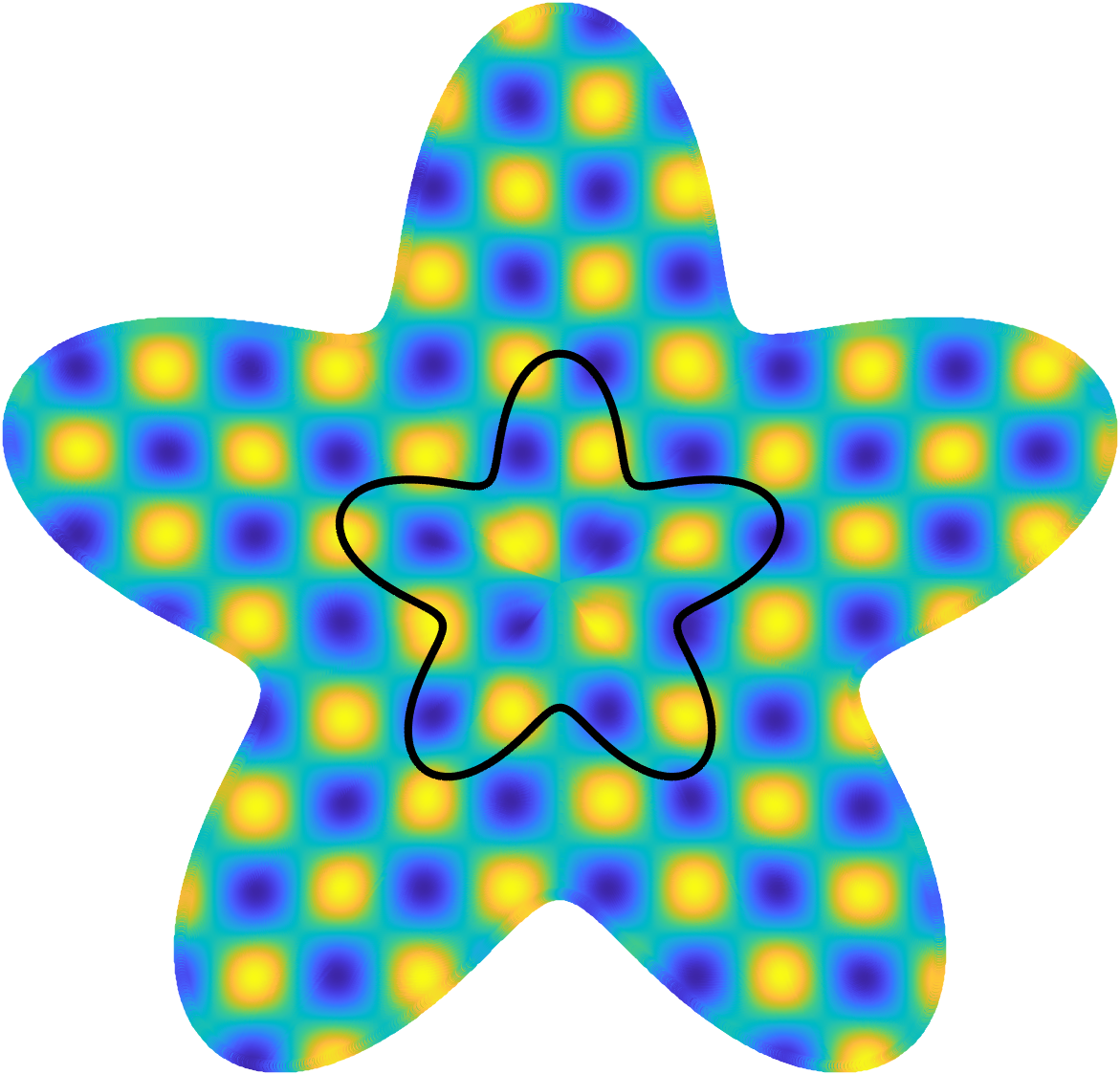}
        \caption{Plot of the extended function obtained by QR-AAA}
        \label{subfig:extendQRAAA}
    \end{subfigure}
    \hfill
    \begin{subfigure}[t]{.32\linewidth}
    \includegraphics[width=\linewidth]{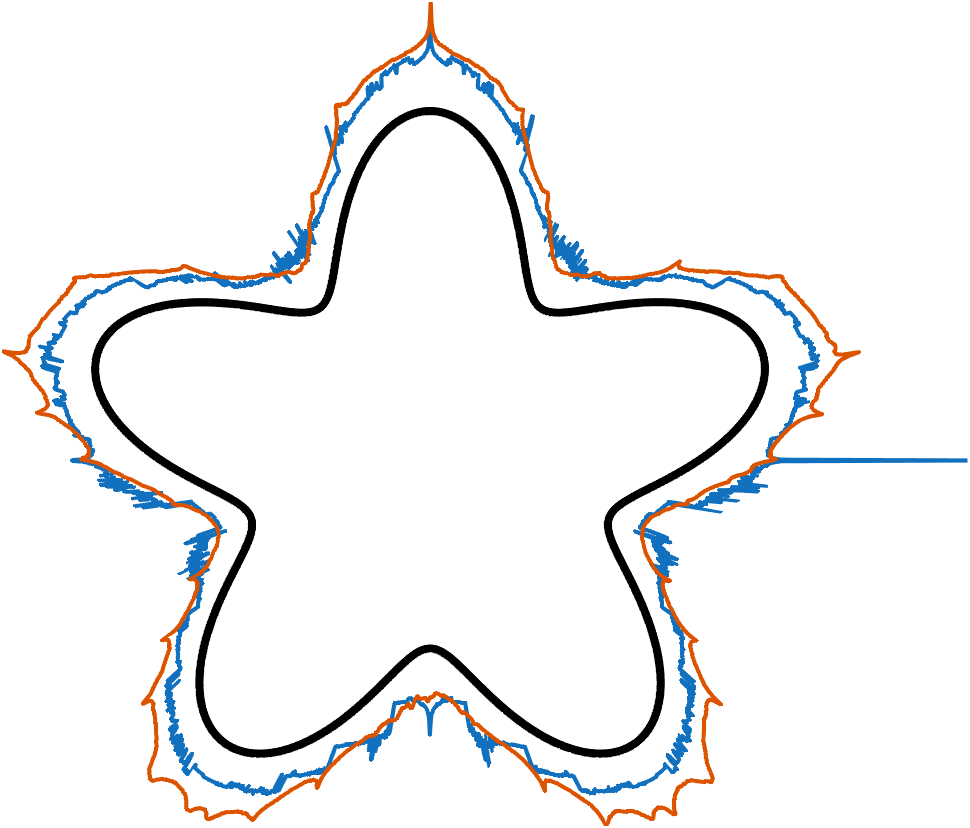}
        \caption{Contour lines showing where ray-AAA (blue) and QR-AAA (red) lose $8$ digits of accuracy}
        \label{subfig:error1e-8}
    \end{subfigure}
    \hfill
    \begin{subfigure}[t]{.32\linewidth}
    \includegraphics[width=\linewidth]{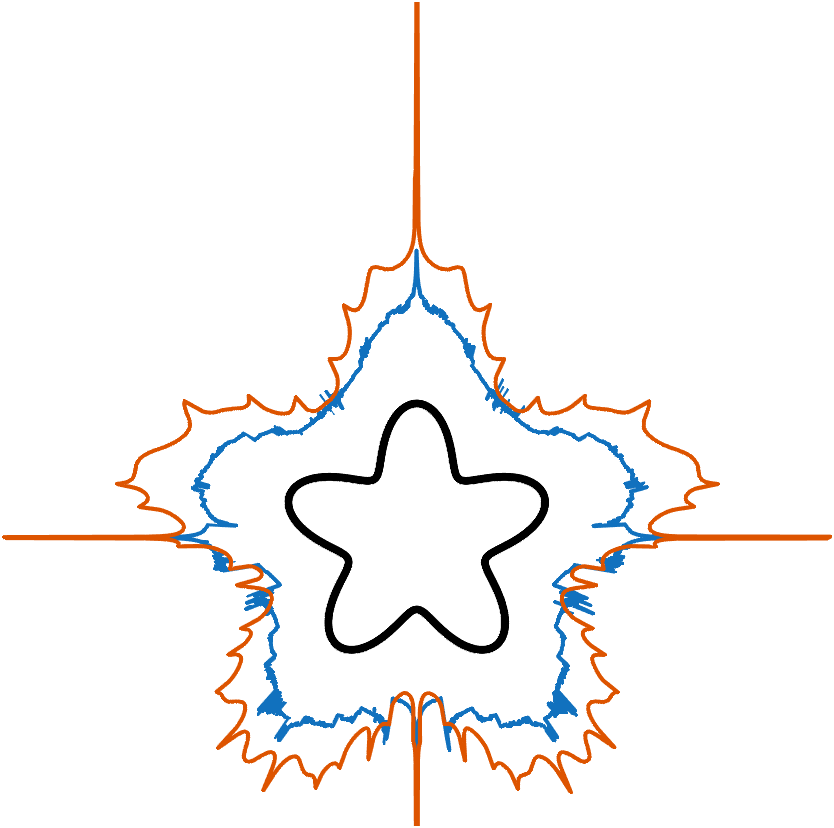}
        \caption{Contour lines showing where ray-AAA (blue) and QR-AAA (red) lose all but $1$ digit of accuracy}
        \label{subfig:extension2}
    \end{subfigure}
    \caption{Plot of the QR-AAA extension of $f(x,y)=-\sin(2\pi x)\sin(2\pi y)$, together with contour lines where single precision and one digit of precision are obtained. The boundary of the original starfish domain is always shown in black.}
    \label{fig:rayqrAAA}
\end{figure}

The next step is natural: model the $2m+1$ functions $\{c_{k}(\rho)\}_{k=-m}^m$ using QR-AAA. In this way, if a successful vector-valued approximation $\mathbf{r}(\rho)=[r_{-m}(\rho),\ldots,r_{m}(\rho)]\approx [c_{-m}(\rho),\ldots,c_{m}(\rho)]$ is obtained, we can extend $f(\rho,\theta)$ beyond $\Omega$ (i.e, $\rho>1$) by setting the extension $\widetilde{f}(\rho,\theta)$ to be
$$\widetilde{f}(\rho,\theta):=r(\rho,\theta) = \sum_{k=-m}^mr_{k}(\rho)e^{i k \theta}.$$ This expression can be evaluated for arbitrary $(\rho,\theta)$, not just for $\theta$ in the original rays of approximation. Since the approximation is bandlimited, it is automatically analytic everywhere in $\mathbb{R}^2$. In principle, especially if $m$ is small (say $\approx50$), the QR step can be omitted. In practice, however, we do not have an a priori estimate for the bandwidth of the angular behavior of $f(\rho,\theta)$. In addition, extension by rational functions requires a high precision in the sampling domain, meaning that underestimating this bandwidth severely impacts the final accuracy of the extension. For this reason, we set $m=200$ and chose $500$ equispaced rays in the angular direction, each of which is discretized using $500$ samples in the radial direction. Subsequently, QR-AAA was used for the rational approximation. The results are shown in Figure~\ref{fig:rayqrAAA}

Per construction, the approximation $\widetilde{f}(x,y)$ is real-analytic and bounded, which, perhaps surprisingly, seems to regularize the problem sufficiently so as to push its loss-of-precision contour lines (See Figure~\ref{fig:rayqrAAA}) \textit{slightly} past those of ray-AAA. Both achieve comparable accuracy close to the boundary though, meaning that the QR-AAA based approach outlined here is mainly favorable for $2$ reasons:
\begin{enumerate}
    \item The extension is analytic and can be evaluated at arbitrary points outside $\Omega$.
    \item The QR-AAA based approximation is faster, and more compact, when compared to ray-AAA.
\end{enumerate}
To be explicit, the rank obtained by the QR decomposition was $12$, and the degree of the vector-valued rational approximation was $18$. Compare this to the maximal degree for ray-AAA, which was $13$, and the average degree, which is $10$. The tolerance for both the truncated QR decomposition and the subsequent SV-AAA approximation were set to $10^{-13}$.

A proper timing comparison of course highly depends on the desired resolution in the angular direction. Let us mention simply that, even with the sub-optimal \texttt{MATLAB} QR-decomposition, and even without using the FFT, the total time to compute $\widetilde{f}(\rho,\theta)$ on a laptop was less than $0.3$ seconds.

\begin{figure}
    \centering
    \begin{subfigure}[t]{.32\linewidth}
    \includegraphics[width=.85\linewidth]{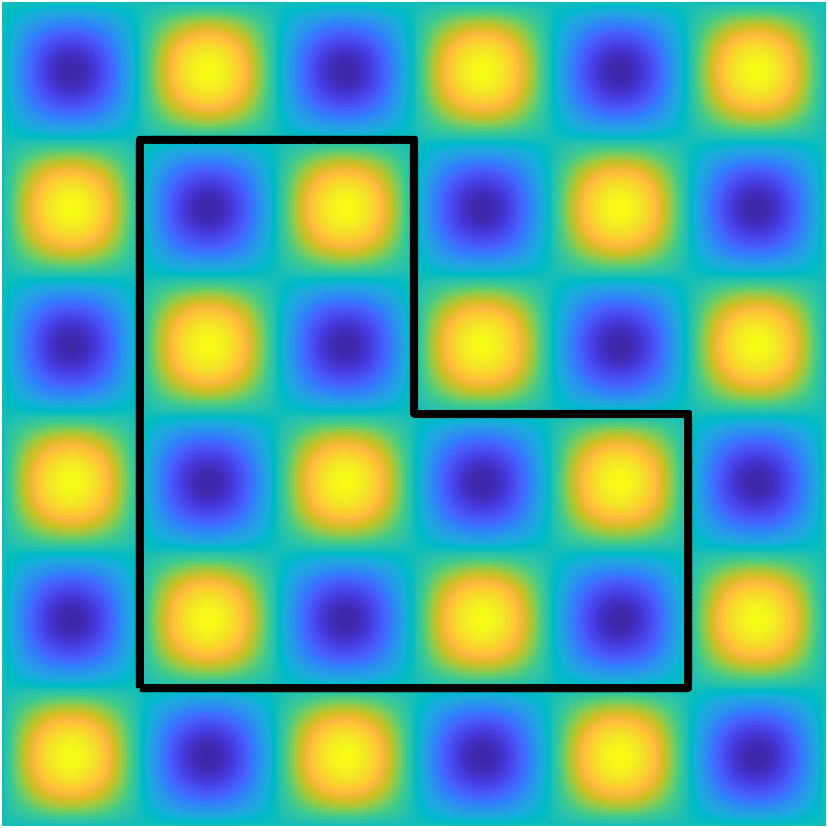}
        \caption{Plot of the function to be extended from an L-shaped domain}
        \label{subfig:extendQRAAAL}
    \end{subfigure}
    \hfill
    \begin{subfigure}[t]{.32\linewidth}
    \centering
    \includegraphics[width=.8\linewidth]{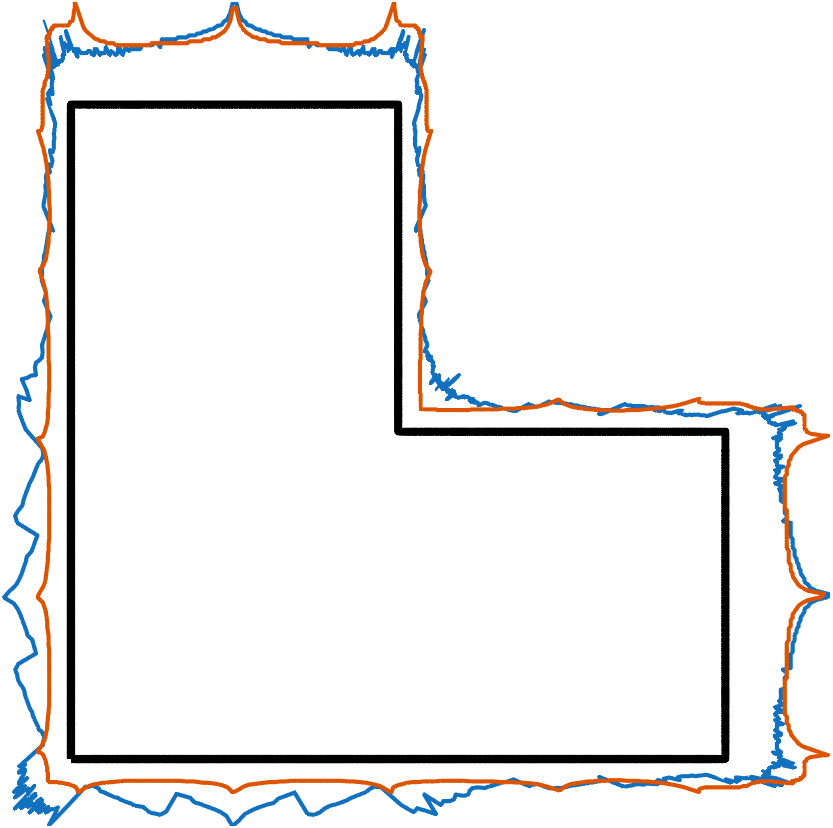}
        \caption{Contour lines showing where ray-AAA (blue) and QR-AAA (red) lose 8 digits of accuracy}
        \label{subfig:error1e-8L}
    \end{subfigure}
    \hfill
    \begin{subfigure}[t]{.32\linewidth}
    \includegraphics[width=\linewidth]{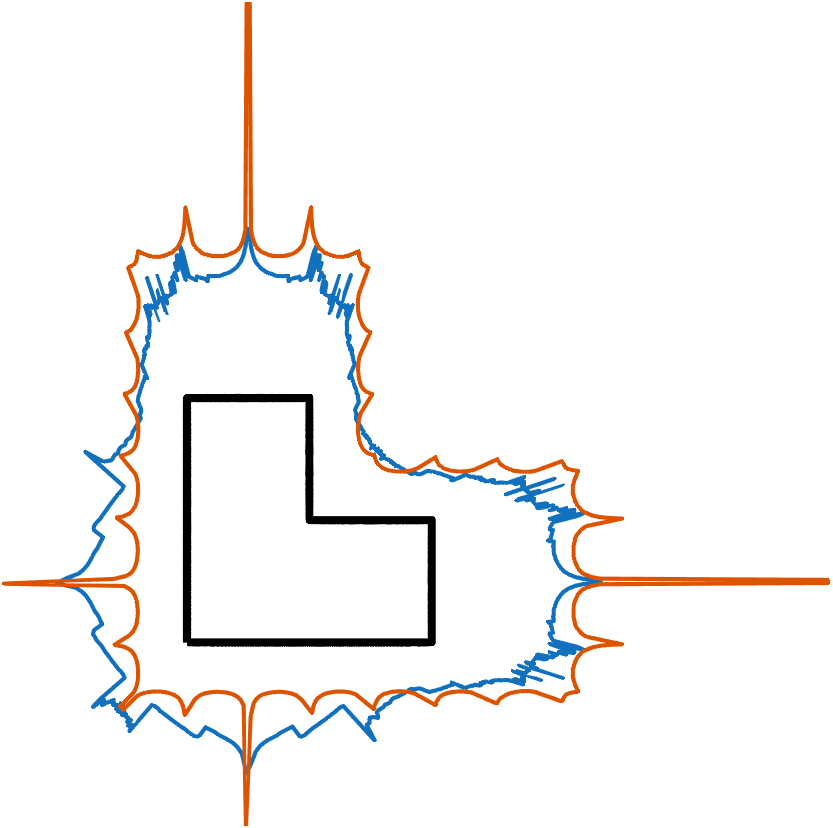}
        \caption{Contour lines showing where ray-AAA (blue) and QR-AAA (red) lose all but 1 digit of accuracy}
        \label{subfig:extensionL}
    \end{subfigure}
    \caption{Plot of the QR-AAA extension of $f(x,y)=-\sin(2\pi x)\sin(2\pi y)$, together with contour lines where single precision and only one digit of precision are obtained. The boundary of the original L-shaped domain is always shown in black.}
    \label{fig:rayqrAAALshape}
\end{figure}

\begin{remark}
    The approach outlined above can be adapted to the case where $\Omega$ is not simply connected, so long as a strip along the boundary inside $\Omega$ can be constructed. We incorporated this in our \texttt{MATLAB} code for the starfish experiment. So long as the width of the strip is sufficiently large, the results are the same and we therefore do not include these here. The case of a non-smooth boundary can also be handled, by switching to, e.g., piecewise polynomial approximations for the angular part of $f(\rho,\theta)$. For an L-shaped domain (which can be found in \texttt{Lshape.m} in our code) the results are shown in Figure~\ref{fig:rayqrAAALshape}. In its most basic form only piecewise differentiability over $\theta$ is guaranteed. Resolving this will be the subject of future work.
\end{remark}

\subsubsection{Extending the norm function}\label{sec:root}
\begin{figure}
    \centering
    \begin{subfigure}[t]{.32\linewidth}
    \includegraphics[width=\linewidth,height=3cm]{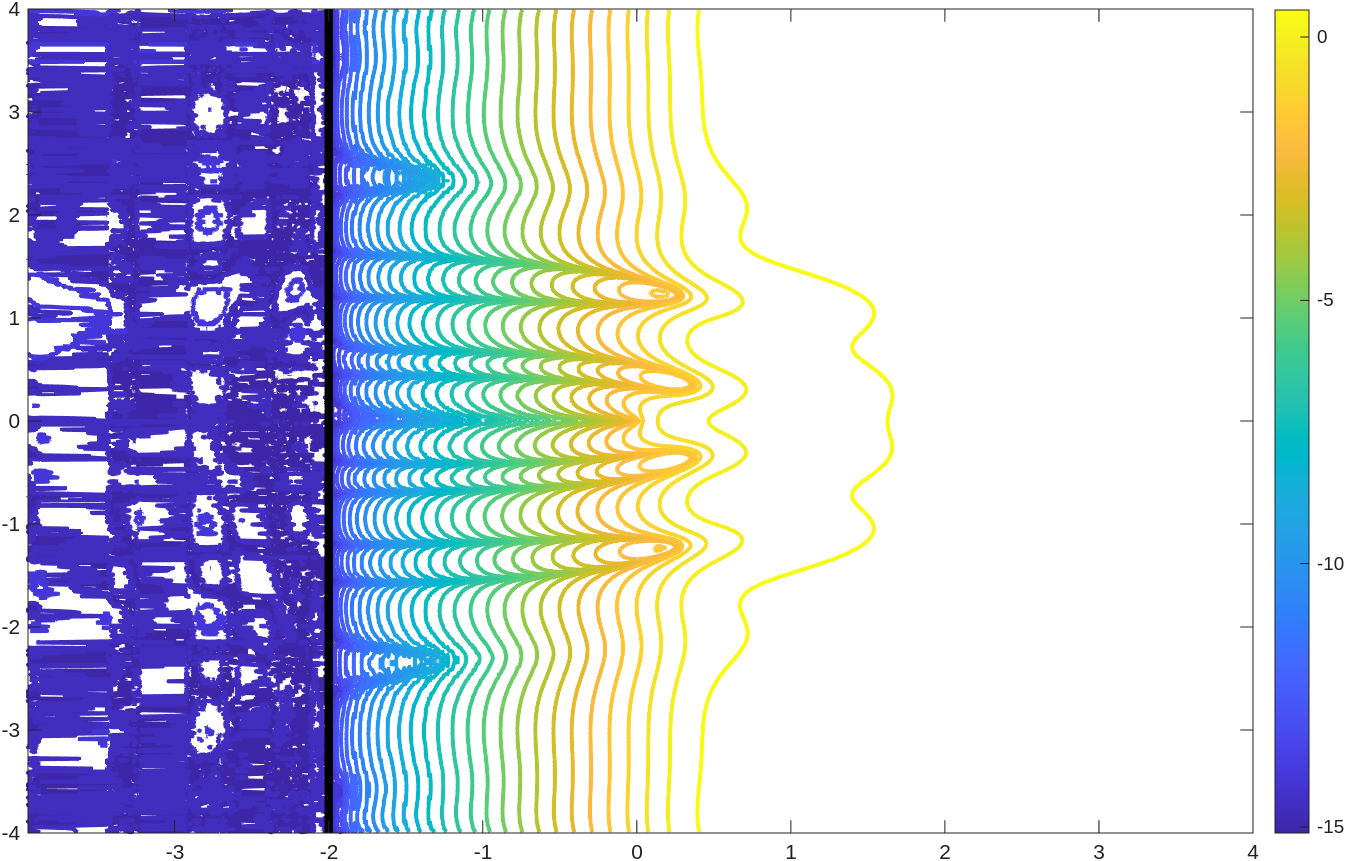}
    \caption{Contour lines for the error of the QR-AAA extension}
    \label{subfig:errQRAAAsqrt}
    \end{subfigure}
    \begin{subfigure}[t]{.32\linewidth}
    \includegraphics[width=\linewidth,height=3cm]{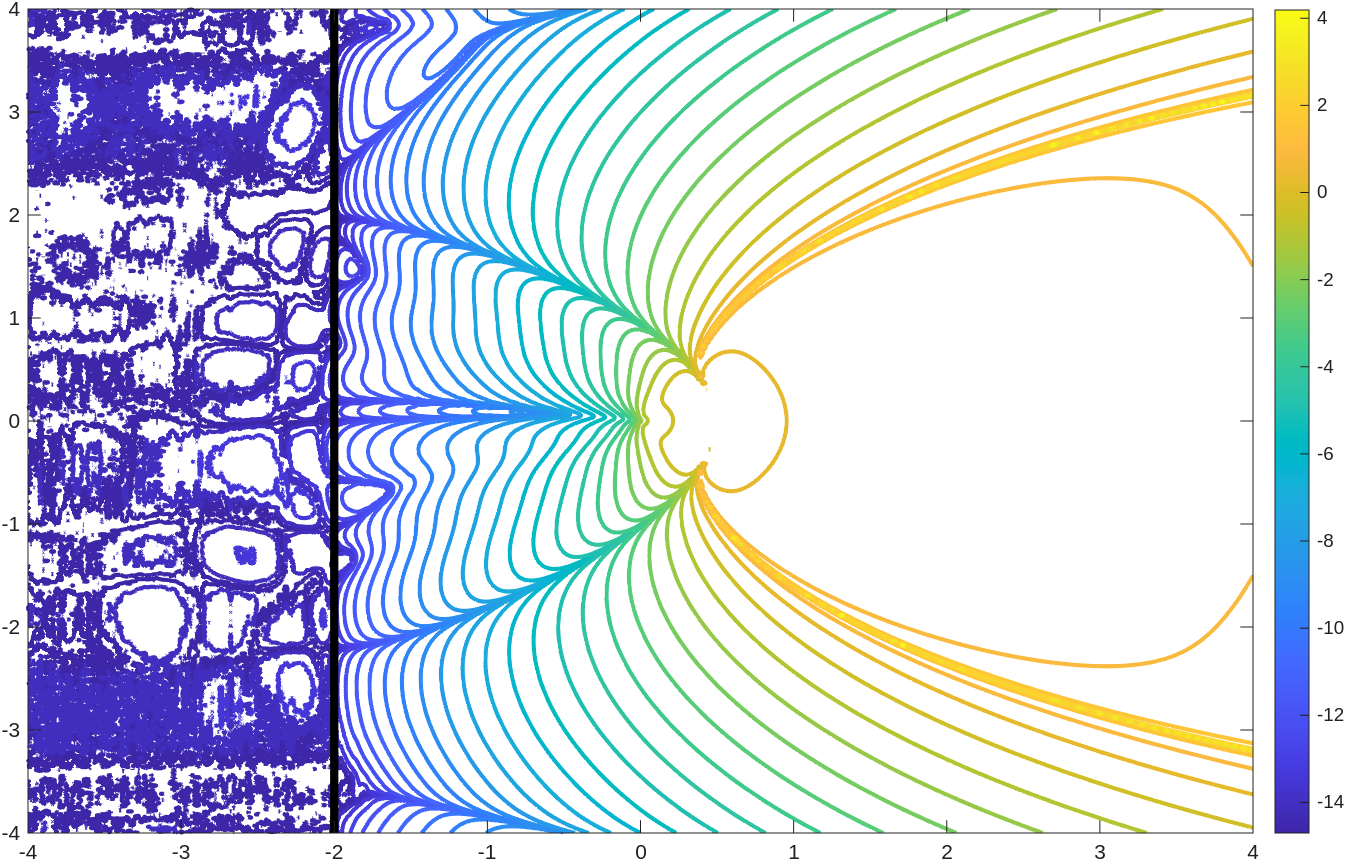}
    \caption{Contour lines for the error of the two-step extension}
    \label{subfig:errpAAAsqrt}
    \end{subfigure}
    \begin{subfigure}[t]{.32\linewidth}
    \includegraphics[width=\linewidth,height=3cm]{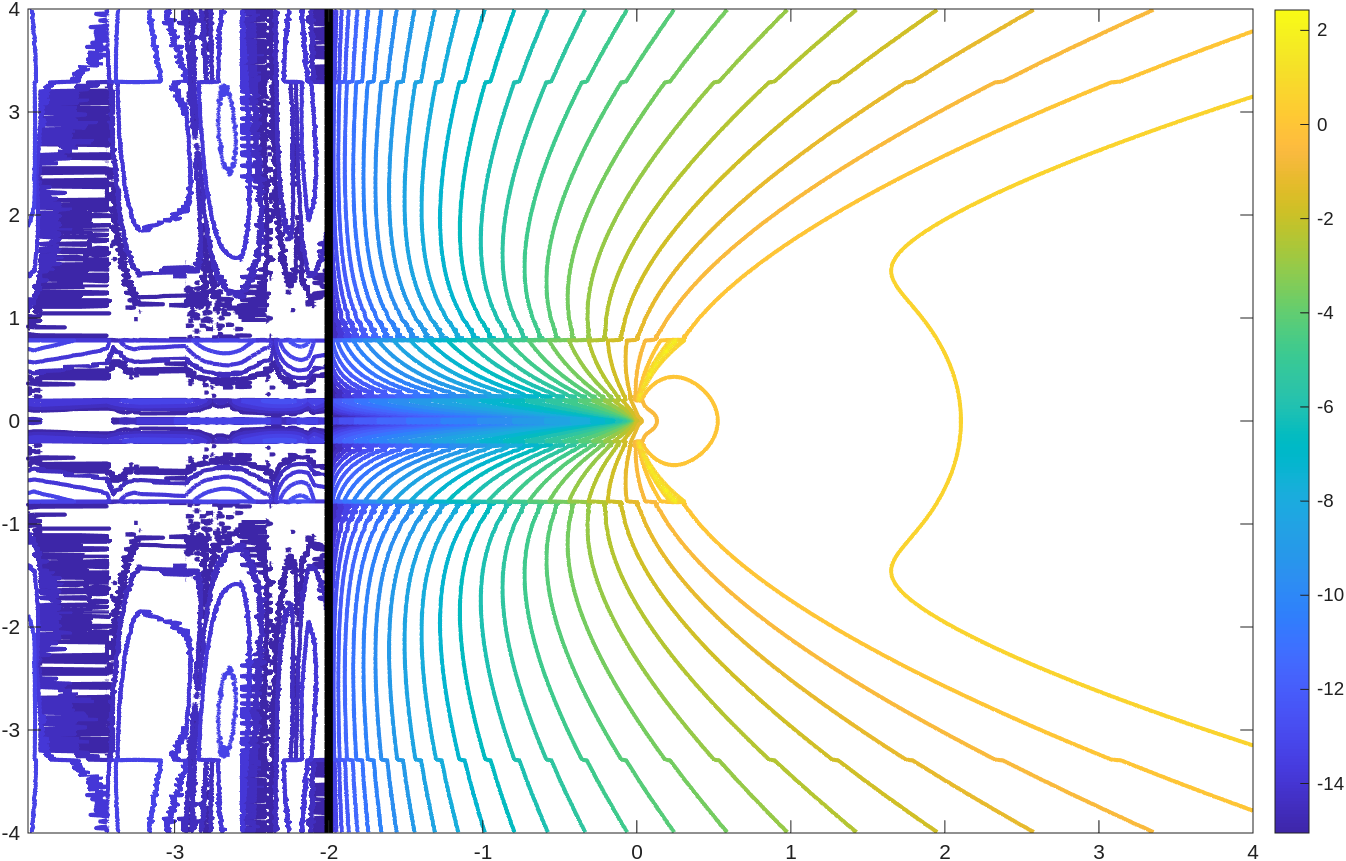}
    \caption{Contour lines for the error of the rayAAA extension}
    \label{subfig:errAAAsqrt}
    \end{subfigure}
    \caption{Errors for the three extension schemes applied to $f(x,y)=\sqrt{x^2+y^2}$.}
    \label{fig:extendSqrt}
\end{figure}

In the previous section, the function $f(x,y)$ could be written as $f(\rho,\theta)$, in such a way that extending $f$ could be achieved by extending $f(\rho,\theta)$ only in the radial direction. In this section we will focus on the function $f(x,y) = \sqrt{x^2+y^2}$, and attempt to extend it in the $x$-direction. Since the underlying function is defined in all of $\R^2$, and its sample domain has a tensor structure, it is more natural to approximate $f(x,y)$ in $\Omega = [-4,-2]\times[-4,4]$ by a \textit{bivariate rational function}. We use the techniques from Section~\ref{sec:multivariate}:
\begin{enumerate}
    \item Separable multivariate rational approximation as introduced in Section~\ref{sec:multivariate}
    \item Two-step p-AAA from Section~\ref{subsubsec:twoStep}.
\end{enumerate}
We set the candidate support point sets $\{Z_x,Z_y\}$ to be finely sampled uniform grids in $[-4,-2]$ and $[-4,4]$ respectively. The results are reported in Figure~\ref{fig:extendSqrt}.

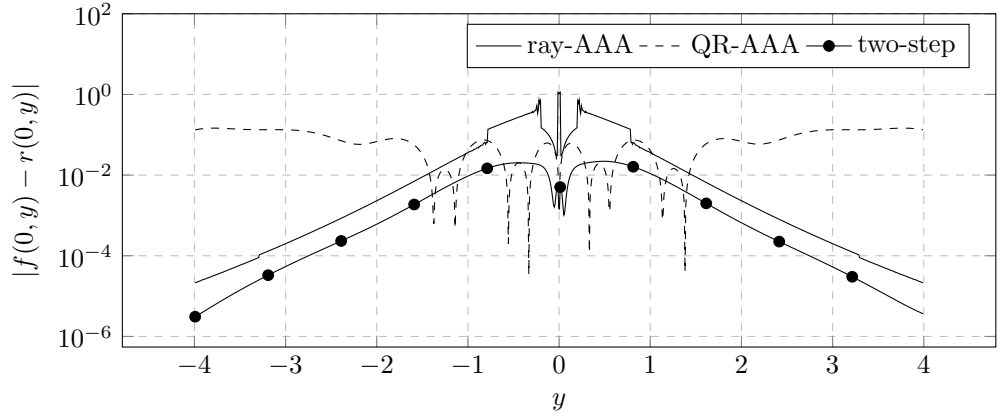
\begin{figure}
    \centering
    \begin{tikzpicture}
    \begin{axis}[
    width=\linewidth,
    height = 6cm,
    ymode=log,
    log basis y=10,
    xlabel={$y$},
    ymax = 100,
    ylabel={$|f(0,y)-r(0,y)|$},
    ymajorgrids=true,
    xmajorgrids=true,
    grid style=dashed,
    legend columns=3,
    legend pos = {north east}
]
\addplot[
    color=black
    ] table [x index = 0 , y index = 1,col sep=comma]{err_sqrt_AAA.tex};
\addlegendentry{ray-AAA}
\addplot[
    color=black,dashed
    ] table [x index = 0 , y index = 1,col sep=comma]{err_sqrt_QRAAA.tex};
\addlegendentry{QR-AAA}
\addplot[
    color=black,mark=*,mark repeat = 100,
    ] table [x index = 0 , y index = 1,col sep=comma]{err_sqrt_pAAA.tex};
    \addlegendentry{two-step}
    \end{axis}
    \end{tikzpicture}
    \caption{Errors of ray-AAA, plain multivariate QR-AAA and two-step pAAA at $x=0$.}
    \label{fig:aroundZero}
    \end{figure}

Because of its separable form, the first approximation performs well only up to \textit{the line} $x=0$, where the singularity is. On the other hand, the more advanced two-step approximation and the ray-AAA approximation suffer similar loss of accuracy around the origin, but maintain a higher accuracy away from the singularity, even at $x=0$. Figure~\ref{fig:aroundZero} shows the error in approximation for $f(0,y)$.

\subsection{Near-field compression using parallel QR-AAA (from \cite{PQRAAA})}\label{sec:BEM}
In this final section, we consider the parallel variant of QR-AAA, \textit{PQR-AAA}. The full details can be found in~\cite{PQRAAA}.

Whenever the number of functions to be approximated is extremely large, the dominant cost in the QR-AAA scheme is the construction of the matrix $\mathbf{F}$ from equation~\eqref{eq:defFmat}. Luckily, parallelizing its construction is often trivial. However, in the distributed setting, a large communication cost would be incurred if one were to straighforwardly apply QR-AAA after this. In fact, it can become completely infeasible to store $\mathbf{F}$ in its entirety. Even in the shared memory setting, it is beneficial to maximally utilize today's parallel computer architecture. Suppose we have $n_{\text{proc}}$ processors available. Briefly, the way PQR-AAA works is:
\begin{enumerate}
    \item Compute $\mathbf{F}$ as $\{\mathbf{F}_{\mu}\}_{\mu=1}^{n_\text{proc}}$ in parallel.
    \item Apply QR-AAA for each processor $\mu$ to obtain $\mathbf{R}_\mu \approx \mathbf{F}_{\mu}$.
    \item Join $\{\mathbf{R}_{\mu}\}_{\mu=1}^{n_{\text{proc}}}$ together into one vector-valued approximation $\mathbf{R}\approx\mathbf{F}$.
\end{enumerate}

As an application, consider the Galerkin boundary element method (BEM) for the Helmholtz equation (see, e.g., \cite{Sauter}, \cite{McLean}). In it, the differential equation in a volume $\Omega$ is transformed into integral equations on the boundary $\partial\Omega$. The simplest boundary integral operator involved, the single layer boundary operator, is defined by
$$(\mathcal{S}(\kappa)[u])(\mathbf{x}) = \lim_{\mathbf{x}\to\partial\Omega}\int_{\partial\Omega}G(\mathbf{x},\mathbf{y};\kappa)u(\mathbf{y})dS_{\mathbf{y}}$$
with $G(\mathbf{x},\mathbf{y};\kappa)$ the Green's kernel at wave number $\kappa$. The operator $\mathcal{S}(\kappa)$ is discretized to $\mathbf{S}(\kappa)$ using Galerkin discretization on a triangular surface mesh. Then $\kappa\mapsto \mathbf{S}(\kappa)$ is a matrix-valued function. We apply PQR-AAA to compress the \textit{near-field}\footnote{The near-field are the degrees-of-freedom that are physically close together.} of $\mathbf{S}(\kappa)$ as it varies over the wave number. This is an important task, as approximating the wavenumber dependence permits one to avoid the costly quadrature needed to construct the near-field of $\mathbf{S}(\kappa)$.

\begin{figure}
\centering
\includegraphics[width = .6\linewidth,height=4cm]{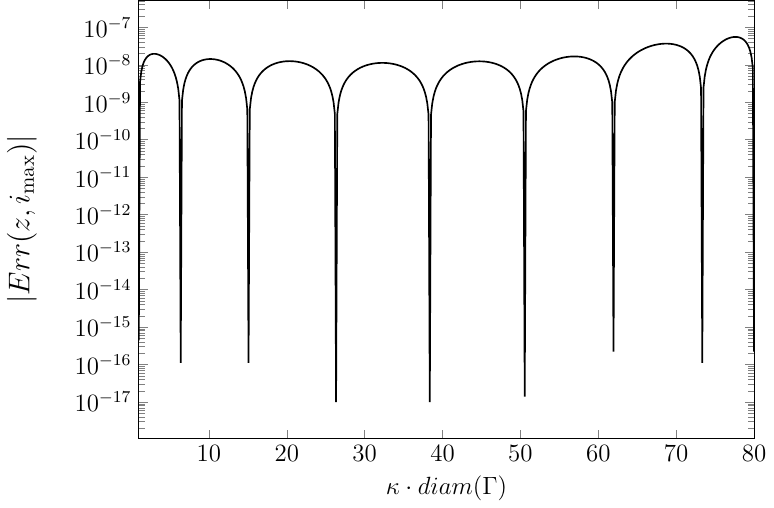}
\caption{Supremum norm error of the PQR-AAA approximation to $S(\kappa)$ over the near-field, varying with the dimensionless wave number $\kappa$}\label{fig:errplotNearField}
\end{figure}

The matrix $\mathbf{F}$ in this case has columns corresponding to $\mathbf{S}(i,j,;Z)$, where $(i,j)$ is an index in the near-field DOFs, and $Z\subset[\kappa_{\min},\kappa_{\max}]$ is a discrete, sufficiently fine subset of the selected frequency range. For our experiment, we discretize $\mathcal{S}$ on a spherical grid obtained from a $6$-fold spherical refinement of an octahedral mesh using the software package \texttt{BEACHpack} (available at \url{https://gitlab.kuleuven.be/numa/software/beachpack}). We ran our experiment over 28 cores. The total number of functions to be approximated is $1000000$, which are distributed (roughly) equally among the cores. The required tolerance was set to $10^{-6}$. The (dimensionless\footnote{The dimensionless wavenumber is the actual wavenumber multiplied by the diameter of $\partial\Omega$}) wavenumber range was set to $[1,80]$, and discretized into $500$ equispaced points. The timings are split up into four constituents:
\begin{enumerate}
\item $t_{\text{F}}$: the maximal time needed to assemble $\mathbf{F}_{\mu}$, over the processors $\mu\in\{1,\ldots,28\}$.
\item $t_{\text{QR}}$: the maximal time needed to compute the approximate pivoted Householder QR decomposition $\mathbf{F}_{\mu}\approx \mathbf{Q}_{\mu}\mathbf{R}_{\mu}$.
\item $t_{\text{AAA}}$: The maximal time needed to compute the SV-AAA approximation of $\mathbf{Q}_{\mu}$.
\item $t_{\text{fin}}$: the time needed to compute the final vector-valued rational approximation for $\mathbf{F}$.
\end{enumerate}
These timings, as well as the supremum norm error over $\mathbf{F}$, are reported for our set-up in table~\ref{tab:timingsQRAAANearfield}. The error over the wavenumber on a fine validation grid is reported in Figure~\ref{fig:errplotNearField}. As can be seen from this figure, the final degree of rational approximation is $12$. The maximum error of the final approximant is well below the requested tolerance.
\begin{table}
\centering
\begin{tabular}{|c|c|c|c|c|}
\hline
$t_{\text{F}}$&$t_{\text{QR}}$&$t_{\text{AAA}}$&$t_{\text{fin}}$&$\mathbf{err}_{\infty}$\\
\hline
$343.213$s&$.926$s&$.019$s&$.259$s&$6.64897\cdot10^{-8}$\\
\hline
\end{tabular}
\caption{Error and timings for the the near-field approximation of $S(\kappa)$.}\label{tab:timingsQRAAANearfield}
\end{table}

\section*{Conclusions}
This manuscript demonstrates that the QR-AAA method is a flexible and effective approach for a wide array of computational applications. For some of the applications outlined in this manuscript, the use of QR-AAA was motivated primarily by its speed and stability. For other applications, new and exciting insights are revealed by the use of vector-valued rational approximation.

\section*{Acknowledgements}
We would like to express our sincere gratitude to Kobe Bruyninckx, for pointing out the subtleties in the transitive property for inexact interpolative decompositions, which helped improve Section~\ref{sec:qr-aaa}. We would also like to thank the referees for their valuable feedback.
\section*{Funding}
The work reported was supported by the National Science Foundation (DMS-2313434) and by the Department of Energy ASCR (DE-SC0025312).
\section*{Data availability}
The methods and implementations in this manuscript are freely available at \url{https://github.com/SimonDirckx/QR-AAA-Applications} and \url{https://github.com/SimonDirckx/pqraaa}. The data for Sections~\ref{sec:Hedgehog}--\ref{sec:extension} can be generated using the scripts in \url{https://github.com/SimonDirckx/QR-AAA-Applications}. The data for Section~\ref{sec:BEM} can be generated using the freely available code at \url{https://gitlab.kuleuven.be/numa/software/beachpack}.
\section*{Conflict of interest}
The authors have no conflict of interest to declare.
\bibliography{main}
\end{document}